\documentclass{article}

\usepackage{arxiv}
\usepackage[numbers,sort,compress]{natbib}
\usepackage[utf8]{inputenc} 
\usepackage[T1]{fontenc}    
\usepackage{hyperref}       
\usepackage{url}            
\usepackage{booktabs}       
\usepackage{amsfonts}       
\usepackage{amssymb}
\usepackage{nicefrac}       
\usepackage{microtype}      
\usepackage{lipsum}		
\usepackage{graphicx}
\usepackage{float}
\usepackage{adjustbox}
\usepackage{subcaption}
\usepackage{mathtools}
\usepackage[capitalize,noabbrev]{cleveref}
\usepackage{multirow}
\usepackage{amsmath}
\usepackage[table,dvipsnames]{xcolor}
\usepackage{amsthm}
\newtheorem{theorem}{Theorem}[section]

\newtheorem{lemma}[theorem]{Lemma}
\newtheorem{corollary}[theorem]{Corollary}

\definecolor{citeblue}{RGB}{0,90,180}
\definecolor{bestgreen}{RGB}{220,240,220}
\definecolor{failred}{RGB}{255,220,220}

\hypersetup{
    colorlinks=true,
    citecolor=citeblue,
    linkcolor=Maroon,
    urlcolor=black
}

\title{Learn the Dynamics, Not the Noise: A Semigroup-Factorized Solver for Stochastic PDEs}

\author{ Bingcheng Hu\thanks{These authors contributed equally.} \\
	Department of Mathematics\\
	Shanghai Normal University\\
	Shanghai, China \\
	\And
	Qi Xu\footnotemark[1] \\
	Department of Mathematics\\
	Shanghai Normal University\\
	Shanghai, China \\
	\AND
	{Zhaoxiang Li }\thanks{Corresponding author,
		email: \texttt{zxli@shnu.edu.cn} }\\
	Department of Mathematics \\
	Shanghai Normal University\\
	Shanghai, China
}

\date{}

\renewcommand{\shorttitle}{Learn the Dynamics, Not the Noise: A Semigroup-Factorized Solver for Stochastic PDEs}

\hypersetup{
pdftitle={A template for the arxiv style},
pdfsubject={q-bio.NC, q-bio.QM},
pdfauthor={David S.~Hippocampus, Elias D.~Striatum},
pdfkeywords={Stochastic Partial Differential Equations, Physics-Informed Learning , Stochastic Convolution ,Semigroup-Based Integration , Neural Time Integration},
}

\begin{document}
\maketitle

\begin{abstract}
Accurate pathwise simulation of stochastic partial differential equations (SPDEs) remains challenging due to the interplay between irregular stochastic forcing and nonlinear dynamics. To address this challenge, we propose \textbf{Bridge-informed Residual Learning in Semigroup Exponential Integration (BRiSE)}, which separates the known stochastic response from the nonlinear dynamics to be learned. BRiSE combines Ornstein--Uhlenbeck bridge-informed path features with a semigroup-factorized decoder, enabling continuous-time trajectory reconstruction with a single network evaluation per interval. The model is trained through pathwise mild-equation residual minimization without requiring reference solution data. Theoretically, we establish decoder consistency and perturbation stability, together with a residual-based global error estimate without an additional stiffness-dependent time-step restriction. A temporal interpolation argument provides a quantitative relationship between the continuous mild loss and the global mean-square trajectory error. Under suitable assumptions, we further establish half-order strong temporal convergence at fixed spectral resolution and a unified error estimate accounting for spectral truncation, stochastic-path approximation, and residual-induced errors. Experiments on six SPDE benchmarks demonstrate that BRiSE consistently achieves the lowest mean space-time relative $L^2$ error among five classical and neural baselines, reducing errors by 11.6\%--89.0\% relative to the strongest competing methods. More broadly, this work suggests a promising paradigm for scientific machine learning: leveraging known mathematical structures to focus neural approximation on the dynamics that genuinely need to be learned.
\end{abstract}


\addtocontents{toc}{\protect\setcounter{tocdepth}{-1}}
\section{Introduction}

\subsection{Background and Motivation}
Stochastic partial differential equations (SPDEs) provide a mathematical framework for modeling spatially distributed systems subject to random perturbations, arising in applications such as fluid dynamics, reaction-diffusion processes, and statistical physics \cite{da2014stochastic,flandoli2011random,cerrai2003stochastic,hairer2014theory}. Unlike deterministic PDEs, SPDEs describe solution trajectories that depend on individual realizations of the driving noise \cite{jentzen2009numerical}. Beyond characterizing statistical quantities such as expectations and variances, accurately resolving these trajectories is essential for understanding how stochastic fluctuations influence the evolution of complex systems. Consequently, accurate pathwise simulation of nonlinear SPDEs constitutes an important computational task in scientific modeling.

Accurate pathwise simulation of nonlinear SPDEs presents substantial mathematical and computational challenges. The temporal irregularity of Wiener forcing makes stochastic solution trajectories difficult to approximate, while nonlinear interactions propagate the effects of random perturbations throughout the evolution \cite{printems2001discretization,kovacs2015backward}. Moreover, spatial differential operators can introduce widely separated time scales, requiring numerical schemes to balance stability and accuracy \cite{abdulle2012numerical,klioba2025pathwise}. These difficulties become particularly pronounced when resolving individual stochastic realizations over extended time horizons. Consequently, the central challenge lies in accurately capturing the interplay between irregular stochastic forcing and nonlinear dynamics while maintaining computational efficiency.

These challenges highlight the need for computational approaches that achieve accurate pathwise simulation without excessive computational cost. In particular, an effective SPDE solver should exploit the mathematical structure of the governing equations while remaining flexible enough to capture nonlinear dynamics under different stochastic realizations. Equally important is the ability to approximate new stochastic trajectories without relying on expensive reference solution data or repeated optimization for each noise realization. Together, these considerations motivate the development of structure-aware and solution-data-free learning frameworks that combine numerical accuracy, computational efficiency, and generalization across stochastic paths.

\subsection{Related Work}

\paragraph{Classical numerical methods for SPDEs.}
Classical numerical methods for SPDEs typically combine spatial discretization with temporal integration schemes, including Euler--Maruyama \cite{kloeden1977numerical}, implicit--explicit methods \cite{ascher1995implicit}, and exponential integrators \cite{jentzen2009overcoming}. While Euler-type methods offer straightforward time stepping, explicit schemes may face restrictive stability constraints when applied to stiff spatial operators \cite{abdulle2022explicit}. Exponential integrators exploit the underlying semigroup structure to treat linear dynamics analytically and incorporate stochastic convolutions. Recent convergence results further demonstrate that sampling stochastic convolutions, rather than relying solely on Wiener increments, can substantially improve temporal convergence rates for certain nonlinear SPDEs with additive noise \cite{djurdjevac2026higher}. These findings underscore the importance of retaining analytically tractable stochastic structures. Nevertheless, accurately approximating the nonlinear drift integral remains challenging, particularly under coarse temporal discretization, where unresolved stochastic fluctuations can significantly influence nonlinear evolution.

\paragraph{Physics-informed and neural SPDE solvers.}
Neural approaches offer flexible alternatives for learning PDE and SPDE solutions. DeepONet \cite{lu2021learning} and FNO \cite{li2021fourier} learn solution operators from data, while PINO \cite{li2024physics} incorporates governing equations into operator learning. For stochastic systems, Neural SPDE \cite{salvi2022neural} learns solution operators conditioned on initial states and noise realizations, accommodating stochastic forcing through a mild-solution-based architecture. More recently, Wiener-chaos-based neural operators \cite{pmlr-v306-shi26i} represent stochastic inputs through chaos features and reconstruct solution trajectories by learning their associated coefficients. Although effective for amortized prediction, these supervised approaches rely on reference solution data, which can be expensive to generate for nonlinear SPDEs. In contrast, physics-informed methods, including PINNs \cite{raissi2019physics}, NN-DO/BO \cite{zhang2020learning}, and NeS-PINN \cite{baptista2026solving}, incorporate governing equations or stochastic modal structures to reduce dependence on labeled solutions. Other approaches address different prediction targets, such as expected solutions in LEC \cite{petursson2025chaos} and probability density evolution in Neural Kolmogorov Equations \cite{bizzi2026neural}. Despite these advances, accurately resolving individual stochastic trajectories without supervised solution data remains a distinct challenge, particularly when nonlinear evolution must be learned consistently with each realization of the driving noise.

\paragraph{Neural time integration and spectral architectures.}
Neural ODEs \cite{chen2018neural} provide a continuous-time modeling framework by parameterizing state dynamics with neural networks and integrating the resulting trajectories numerically. Building on this paradigm, NeuSA \cite{bizzi2026neuro} combines spectral representations with Neural ODEs to learn the evolution of PDE solutions in a finite-dimensional coefficient space. NEXT \cite{marques2026next} further incorporates exponential time differencing, analytically handling stiff linear dynamics while learning the nonlinear remainder. These advances demonstrate the benefits of combining neural approximation with structured numerical integration. However, their formulations primarily target deterministic PDEs, whereas pathwise SPDE simulation introduces an additional challenge: stochastic forcing continuously perturbs the state and influences the nonlinear evolution within each time interval. Consequently, extending neural time integration to SPDEs requires more than incorporating random inputs into learned dynamics; it calls for a principled treatment of stochastic propagation and its interaction with nonlinear drift. This motivates a framework that explicitly preserves analytically tractable stochastic structures while learning only the nonlinear corrections informed by the realized noise.

\begin{figure*}[t]
\centering
\includegraphics[width=\textwidth]{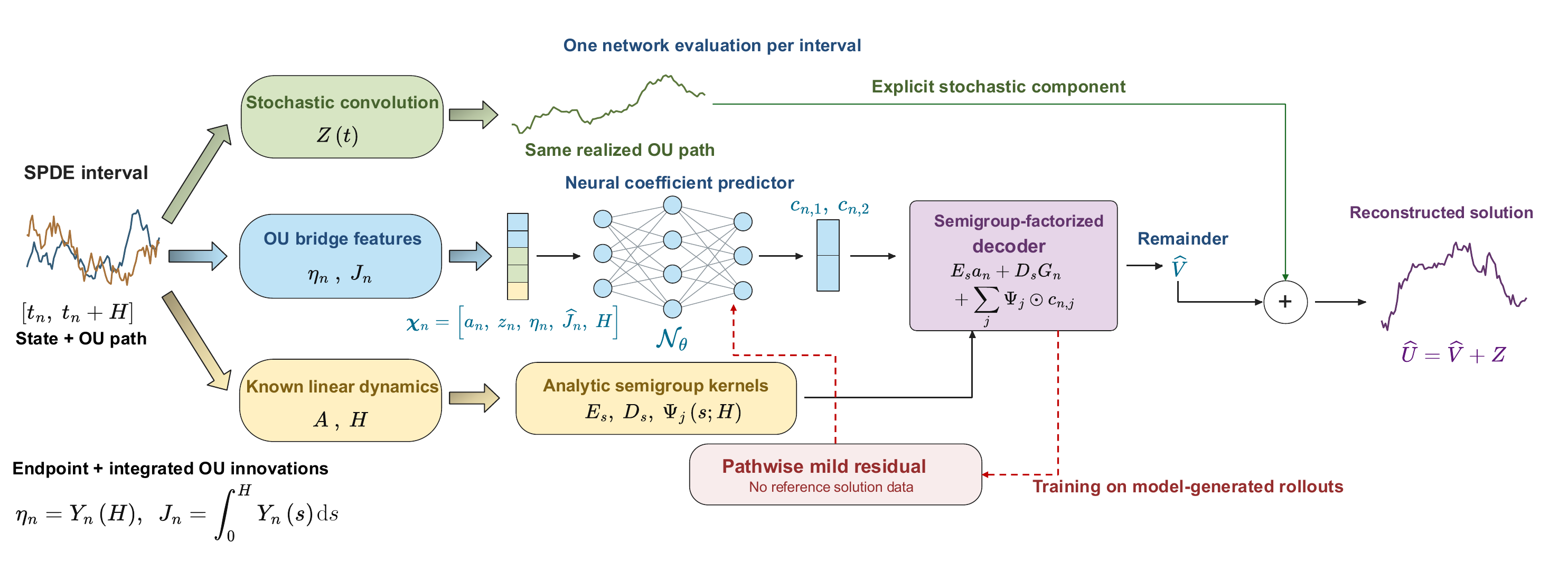}
\caption{BRiSE framework.}
\label{fig:teaser}
\end{figure*}

\subsection{Main Contributions}
To address these challenges, we propose  \textbf{Bridge-informed Residual Learning in Semigroup Exponential Integration (BRiSE)}, a structure-aware neural solver for nonlinear SPDEs with additive noise that combines stochastic path encoding, semigroup-factorized time integration, and physics-informed learning, whose framework is shown in Figure~\ref{fig:teaser}. The key idea is to separate what can be computed analytically from what must be learned. Rather than learning stochastic solution trajectories from expensive reference data, BRiSE explicitly accounts for the stochastic convolution and focuses neural approximation on path-dependent nonlinear dynamics. By combining the analytical structure of exponential integration with learnable temporal corrections, BRiSE improves the approximation of nonlinear evolution beyond prescribed time-stepping updates. Moreover, its physics-informed training requires no reference solution trajectories and enables inference across unseen noise realizations without path-specific retraining. Our main contributions are summarized as follows:
\begin{itemize}
\item We introduce an Ornstein--Uhlenbeck (OU) bridge-informed encoding strategy that combines endpoint and integrated innovations to capture complementary information about within-interval stochastic forcing, enabling path-conditioned neural time integration.
\item We develop a semigroup-factorized neural decoder that parameterizes nonlinear temporal corrections through path-conditioned coefficients and analytically integrates them against the known semigroup, enabling continuous-time trajectory reconstruction with a single network evaluation per interval.
\item We establish a solution-data-free learning framework that trains the proposed decoder by enforcing pathwise mild-equation consistency along model-generated rollouts, using the same stochastic realizations for path encoding and residual evaluation without requiring reference solution trajectories.
\item We establish a theoretical framework that separates spectral truncation, stochastic-path approximation, and residual-induced errors, together with stability and conditional convergence results that characterize how mild-equation residuals control global pathwise accuracy without stiffness-dependent time-step restrictions.
\end{itemize}

\subsection{Organization}
This paper is organized as follows: Section~\ref{sec:preliminaries} introduces the mathematical preliminaries, followed by the proposed BRiSE framework in Section~\ref{sec:method}. Section~\ref{sec:theory} develops the theoretical analysis and proves the guarantees stated in Section~\ref{sec:method}. Section~\ref{sec:numerical_exp} presents numerical experiments, while Section~\ref{sec:ab} provides hyperparameter, ablation, and noise robustness studies. Section~\ref{sec:conclusion} closes the paper, with additional details provided in the Appendix.

\section{Preliminaries}\label{sec:preliminaries}
In this section, we introduce the mathematical foundations of BRiSE, including the mild formulation of semilinear SPDEs, stochastic convolution decomposition, and spectral Galerkin representation.

\subsection{Problem Setup and Stochastic Convolution}
We consider a class of semilinear SPDEs with additive noise, defined on a Hilbert space $\mathcal H=L^2(D)$, where $D\subset\mathbb R^{d_{\mathrm{sp}}}$ is a bounded spatial domain and $d_{\mathrm{sp}}$ is the spatial dimension. The governing equation is given by
\begin{equation}\label{eq:spde}
\begin{aligned}
&\mathrm{d}u(t)=\left(Au(t)+F(u(t))\right)\,\mathrm{d}t+\sigma\,\mathrm{d}W_t^Q, \\
&u(0)=u_0,
\end{aligned}
\end{equation}
where $t\in[0,T]$, $A$ is a linear operator generating a strongly continuous semigroup $S(t)=e^{tA}$, $F$ denotes the nonlinear drift, and $W_t^Q$ is a $Q$-Wiener process with covariance operator $Q$. The parameter $\sigma$ controls the noise intensity. We assume prescribed homogeneous or periodic boundary conditions and that the SPDE is well-posed in the considered spatial setting.

Applying the variation-of-constants formula \cite{da2014stochastic} to Equation~\eqref{eq:spde} gives the mild formulation
\begin{equation}\label{eq:mild_spde}
u(t)= S(t)u_0 + \int_0^t S(t-s)F(u(s))\,\mathrm{d}s +\sigma\int_0^t S(t-s)\,\mathrm{d}W_s^Q.
\end{equation}
This representation decomposes the solution into three components, which are the linear propagation of the initial condition, the accumulated nonlinear drift, and the stochastic convolution induced by the noise. In particular, the semigroup $S(t)$ governs the propagation of both deterministic and stochastic contributions, providing a natural analytical framework for separating the stochastic forcing from the nonlinear dynamics.

To decouple the stochastic forcing from the nonlinear dynamics, we introduce the stochastic convolution $Z(t)$ and the remainder $v(t)$, defined as
\begin{equation}\label{eq:split}
Z(t)=\sigma\int_0^t S(t-s)\,\mathrm{d}W_s^Q, \quad v(t)=u(t)-Z(t).
\end{equation}
Substituting $u=v+Z$ from Equation~\eqref{eq:split} into Equation~\eqref{eq:mild_spde} and subtracting $Z(t)$ gives
\begin{equation}\label{eq:mild_remainder}
v(t)=S(t)u_0+ \int_0^t S(t-s)F\big(v(s)+Z(s)\big)\,\mathrm{d}s.
\end{equation}
This decomposition isolates the stochastic convolution, which is determined by the known linear operator and noise process, leaving a nonlinear evolution driven by the realized stochastic path. Consequently, the remainder can be treated pathwise without differentiating the irregular noise trajectory. This representation provides the mathematical foundation for the noise-conditioned neural time integration developed in BRiSE.

\subsection{Spectral Galerkin Representation}
To obtain a finite-dimensional representation of the remainder dynamics, we employ a spectral Galerkin approximation. We assume that the linear operator $A$ and the covariance operator $Q$ admit a common orthonormal eigenbasis $\{\phi_k\}_{k\geq 1}$ of $\mathcal H$ \cite{brehier2018weak}, satisfying
\begin{equation}\label{eq:eigenbasis}
A\phi_k=-\lambda_k\phi_k, \quad Q\phi_k=q_k\phi_k,
\end{equation}
where $\lambda_k\geq 0$ and $q_k\geq 0$. Truncating the spectral expansions to $K$ modes, we represent the remainder and stochastic convolution as
\begin{equation}\label{eq:galerkin_expansions}
v_K(t)=\sum_{k=1}^{K}a_k(t)\phi_k, \quad Z_K(t)=\sum_{k=1}^{K}z_k(t)\phi_k.
\end{equation}
Substituting the truncated expansions from Equation~\eqref{eq:galerkin_expansions} into Equation~\eqref{eq:mild_remainder} and projecting onto $\phi_k$ gives a scalar integral equation with kernel $e^{-\lambda_k(t-s)}$ by Equation~\eqref{eq:eigenbasis}. Differentiating in time and collecting the $K$ coefficients gives
\begin{equation}\label{eq:galerkin_ode}
\dot a(t)=-\Lambda a(t)+G(a(t),z(t)),
\end{equation}
where $a,z\in\mathbb R^K$, $\Lambda=\operatorname{diag}(\lambda_1,\ldots,\lambda_K)$, and the nonlinear term is defined componentwise by
\begin{equation}
G_k(a,z)=\left\langle F\bigg(\sum_{j=1}^{K}(a_j+z_j)\phi_j\bigg), \phi_k \right\rangle_{\mathcal H},
\end{equation}
where $\left\langle \cdot,\cdot \right\rangle_{\mathcal H}$ denotes the inner product in $\mathcal H$.
Throughout, $\|u\|$ denotes the spatial $L^2$ norm of $u\in\mathcal H$, and $\|x\|_2$ denotes the Euclidean norm of a finite-dimensional vector $x$. For a matrix $B$, $\|B\|_{\mathrm{op}}$ is the induced Euclidean operator norm. For paths on a time interval $I$, we write
\[
  \|w\|_{\infty,I}=\sup_{t\in I}\|w(t)\|_2
  \quad\text{for coefficient paths},\qquad
  \|u\|_{\infty,I;\mathcal H}=\sup_{t\in I}\|u(t)\|
  \quad\text{for spatial fields}.
\]
For a random vector $X$, we write $\|X\|_{L^2(\Omega)}=(\mathbb E[\|X\|_2^2])^{1/2}$, using $|X|$ in place of $\|X\|_2$ for a random scalar. In particular, the root-mean-square error of a coefficient path $w$ on $I$ is $(\mathbb E[\|w\|_{\infty,I}^2])^{1/2}$.
This formulation reduces the infinite-dimensional SPDE remainder dynamics to a finite-dimensional stochastically driven system, providing the foundation for BRiSE to construct a neural time integrator in the spectral coefficient space.

Substituting the stochastic convolution in Equation \eqref{eq:galerkin_expansions} into Equation \eqref{eq:split} and  projecting the resulting expression onto $\phi_k$, we obtain $z_k(t)=\sigma\sqrt{q_k}\int_0^t e^{-\lambda_k(t-s)}\,d\beta_k(s)$. Its differential form is
\begin{equation}\label{eq:modal_ou}
\mathrm{d}z_k(t)=-\lambda_k z_k(t)\,\mathrm{d}t+\sigma\sqrt{q_k}\,d\beta_k(t),\quad z_k(0)=0,
\end{equation}
where $\{\beta_k\}_{k=1}^{K}$ are independent standard Brownian motions. This representation characterizes the stochastic forcing in the spectral domain and provides the basis for the interval-wise OU path encoding introduced in Section~\ref{subsec:OU Bridge-Informed Path Encoding}. Building on this spectral formulation, we next introduce BRiSE, which leverages the OU path structure and known semigroup dynamics to construct a neural time integrator.

\section{Methodology}\label{sec:method}
In this section, we propose BRiSE, which combines OU bridge-informed path features, a semigroup-factorized decoder, and pathwise residual training. We describe each component and discuss its implications for numerical accuracy. The corresponding analytical framework and theoretical proofs are presented in Section~\ref{sec:theory}.

\subsection{OU Bridge-Informed Path Encoding} \label{subsec:OU Bridge-Informed Path Encoding}
To characterize the stochastic forcing within each integration interval, we partition $[0,T]$ into $N$ coarse intervals $[t_n,t_{n+1}]$, with length $H=T/N$ and nodes $t_n=nH$. Let $z_n=z(t_n)$ denote the stochastic convolution coefficients at the beginning of the interval, and define $E_s=e^{-\Lambda s}$. Solving Equation~\eqref{eq:modal_ou} from $t_n$ to $t_n+s$ gives the propagated state $E_sz_n$ plus a stochastic integral over this interval. We define the latter as the OU innovation
\begin{equation}\label{eq:innovation}
Y_n(s)=z(t_n+s)-E_sz_n,\quad s\in[0,H].
\end{equation}
By subtracting the deterministic propagation $E_sz_n$, this representation isolates the stochastic forcing accumulated within the current interval.

To capture both the terminal and accumulated effects of the stochastic forcing, we take the endpoint and the time integral of the innovation in Equation~\eqref{eq:innovation} with
\begin{equation}\label{eq:ou_features}
\eta_n=Y_n(H),\quad J_n=\int_0^H Y_n(s)\,\mathrm{d}s.
\end{equation}
Given $Y_n(H)=\eta_n$, the conditional law of $Y_n$ is an OU bridge from zero to $\eta_n$. The integral $J_n$ contains the integral of the conditional mean and the integral of the centered bridge fluctuations. Thus $\eta_n$ describes the terminal forcing, while $J_n$ also captures fluctuations within the interval. BRiSE uses the integral of the realized innovation directly as its bridge-informed feature.

Let $x_n$ denote the predicted remainder coefficients at $t_n$, with $x_0=a(0)$ for exact initialization. We combine this state with the stochastic features and interval length to form the neural input
\begin{equation}
\boldsymbol{\chi}_n=[x_n,z_n,\eta_n,\widehat J_n,H] \in \mathbb R^{4K+1},
\end{equation}
where $\widehat J_n$ denotes the numerical approximation of $J_n$. All stochastic features are extracted from the same OU realization used in the pathwise residual evaluation. This representation enables BRiSE to condition its nonlinear temporal corrections on both the current system state and the stochastic forcing within each interval. Further details on OU path sampling and the numerical construction of the integrated innovation are provided in Appendix~\ref{app:OU Path Construction}. Building on this path-conditioned representation, we next introduce a semigroup-factorized neural decoder that reconstructs continuous-time trajectories from interval-wise correction coefficients.

\subsection{Semigroup-Factorized Neural Decoder}\label{sec:brise_theory}
Building on the OU path-conditioned representation, we develop a semigroup-factorized neural decoder that reconstructs continuous-time trajectories through a single coefficient prediction per integration interval. Applying the variation-of-constants formula to Equation~\eqref{eq:galerkin_ode} on $[t_n,t_n+s]$ gives
\begin{equation}\label{eq:interval_mild}
a(t_n+s) = E_sa(t_n)+ \int_0^s E_{s-r} G\big(a(t_n+r),z(t_n+r)\big)\,\mathrm{d}r.
\end{equation}
For an interval starting from the predicted state $x_n$, BRiSE approximates the nonlinear forcing by a polynomial in the normalized time $r/H$ by
\begin{equation}\label{eq:forcing_poly}
G_n+\sum_{j=1}^{M}p_j(r/H)c_{\theta,n,j}.
\end{equation}
Here $G_n=G(x_n,z_n)$, $p_j(\tau)=\tau^j$, and $M$ is the highest correction degree. The neural network $N_\theta$, with trainable parameters $\theta$, maps $\boldsymbol{\chi}_n\in\mathbb R^{4K+1}$ to the tuple of coefficient vectors $(c_{\theta,n,j})_{j=1}^{M}\in(\mathbb R^K)^M$. Substituting Equation~\eqref{eq:forcing_poly} for the nonlinear forcing in Equation~\eqref{eq:interval_mild}, using $x_n$ as the initial state, and integrating term by term gives the decoder
\begin{equation}\label{eq:brise_decoder}
	\widehat a_{\theta,n}(s) = E_sx_n+D_sG_n+ \sum_{j=1}^{M} \Psi_j(s;H)\odot c_{\theta,n,j},
\end{equation}
where the constant and polynomial forcing terms have respective weights
\begin{equation}\label{eq:basis_integrals}
D_s=\int_0^s E_{s-r}\,\mathrm{d}r,\quad \Psi_{j,k}(s;H)=\int_0^s e^{-\lambda_k(s-r)}p_j(r/H)\,\mathrm{d}r.
\end{equation}
Here $\Psi_j(s;H)=(\Psi_{j,k}(s;H))_{k=1}^K$ and $\odot$ denotes elementwise multiplication. We adopt $M=2$, corresponding to linear and quadratic temporal corrections. This factorization incorporates the known semigroup dynamics while allowing the neural coefficients to adapt to individual noise realizations. The resulting decoder reconstructs continuous-time trajectories with a single network evaluation per interval, without additional neural evaluations at intermediate times. Further details on the semigroup-integrated temporal basis functions are provided in Appendix~\ref{app:temporal_basis}. We next introduce a pathwise mild-equation learning objective to train these correction coefficients without reference solution trajectories.

The decoder preserves the prescribed interval state $x_n$ and its initial drift $-\Lambda x_n+G_n$. Under the trajectory-wise Lipschitz bounds in Section~\ref{sec:theory_setting}, its sensitivity to the initial state, noise, and correction coefficients is bounded by Lemma~\ref{lem:theory_decoder}. The semigroup basis constants introduce no dependence on the largest eigenvalue. The complete statement and proof are in Section~\ref{sec:theory_decoder}.

\subsection{Solution-Data-Free Pathwise Learning}\label{sec:pathwise_learning}
To train BRiSE without reference solution trajectories, we introduce a physics-informed learning objective that enforces the pathwise mild formulation on the reconstructed trajectories. We insert the decoder from Equation~\eqref{eq:brise_decoder} into Equation~\eqref{eq:interval_mild} with initial state $x_n$ and subtract the right-hand side from the left to define the residual
\begin{equation}\label{eq:brise_residual}
\mathcal R_{\theta,n}(s)={}\widehat a_{\theta,n}(s)-E_s x_n -\int_0^s E_{s-r}G\big(\widehat a_{\theta,n}(r),z(t_n+r)\big)\,\mathrm{d}r.
\end{equation}
To relate the residual to rollout error, let $z$ be the exact stochastic convolution path in coefficient coordinates and $\widetilde z$ the path supplied to the numerical method. For an inexact path, $\mathcal R_{\theta,n}$ is evaluated with $\widetilde z$ in place of $z$, and the decoder and its inputs use that same path. Let $L_a,L_z$ be the Lipschitz bounds for $G$ in its remainder and noise arguments, respectively, as specified in Section~\ref{sec:theory_setting}.
\begin{theorem}[Pathwise residual-based error estimate]\label{thm:theory_residual}
In the setting of Section~\ref{sec:theory_setting}, let $a$ be the exact Galerkin remainder and $\widehat a$ the rollout driven by $\widetilde z$. Set $e_0=a(0)-\widehat a(0)$ and $\rho_n=\|\mathcal R_{\theta,n}\|_{\infty,[0,H]}$ with the residual evaluated along $\widetilde z$.
Then
\begin{equation}
\label{eq:theory_pathwise_error}
\|a-\widehat a\|_{\infty,[0,T]} \le e^{L_aT}\Bigl(\|e_0\|_2+\sum_{n=0}^{N-1}\rho_n+L_z\int_0^T\|z(t)-\widetilde z(t)\|_2\,\mathrm{d}t\Bigr).
\end{equation}
\end{theorem}
Subtracting the residual identity from Equation~\eqref{eq:interval_mild} and applying Gr\"onwall's inequality over successive intervals gives the bound above. The proof and mean-square extension are given in Section~\ref{sec:theory_decoder}.

Although the decoder analytically incorporates the semigroup structure, the residual is generally nonzero because its learned temporal forcing approximation need not coincide with the physical nonlinearity. We therefore minimize this residual to learn correction coefficients that are consistent with the governing dynamics.

Applying quadrature to the integral in Equation~\eqref{eq:brise_residual} gives the numerical residual $\widehat{\mathcal R}_{\theta,n}(s)$. Given $B$ sampled noise paths, $N$ coarse intervals, and $N_{\mathrm c}$ collocation points per interval, the training objective is defined as
\begin{equation}
\begin{aligned}
\mathcal L_{\mathrm{mild}}(\theta) &= \mathbb E_{W^Q} \Bigg[ \frac{1}{N}\sum_{n=0}^{N-1} \mathbb E_{s\sim\mathcal U(0,H)} \bigg[ \frac{\|\mathcal R_{\theta,n}(s)\|_2^2}{H^2} \bigg] \Bigg]
\approx \frac{1}{B N N_{\mathrm c}} \sum_{b=1}^{B} \sum_{n=0}^{N-1} \sum_{m=1}^{N_{\mathrm c}} \frac{ \left\| \widehat{\mathcal R}_{\theta,n}^{(b)}(s_m) \right\|_2^2 }{H^2}.
\end{aligned} \label{eq:mild_loss}
\end{equation}
Here $\mathbb E [\cdot ]$ denotes expectation. The outer expectation is taken over realizations of the driving Wiener process $W^Q$, and the inner expectation is taken over the local time $s$ sampled uniformly from $(0,H)$ for each fixed path. In the discrete loss, $s_m$ is the $m$-th collocation time and the superscript $(b)$ identifies the $b$-th sampled path. Dividing the residual by $H$ before squaring measures the defect per unit interval length.

Crucially, the stochastic path used for residual evaluation is identical to that used for OU feature extraction, ensuring consistent pathwise constraints. During training, the predicted endpoint initializes the next interval via $x_{n+1}=\widehat a_{\theta,n}(H)$. This allows BRiSE to optimize its temporal corrections along its own predicted trajectories, without relying on reference solutions for training or checkpoint selection. Further details on residual quadrature and rollout-based optimization are provided in Appendix~\ref{app:numerical_residual_evaluation}.

Once trained, BRiSE advances the solution through successive coarse intervals using the predicted endpoint states, while reconstructing the full stochastic trajectory by combining the neural remainder with the corresponding stochastic convolution.

The continuous loss also controls rollout error when residual derivatives are bounded. Write $J=\mathcal L_{\mathrm{mild}}(\theta)$ and define
$$ \overline{\mathcal D}_R=\left(\mathbb E\bigg[\max_n\sup_{s\in[0,H]} \|\mathcal R_{\theta,n}'(s)\|_2^2\bigg]\right)^{1/2}. $$
The next theorem follows by using the residual derivative bound to estimate the sum of interval suprema in Equation~\eqref{eq:theory_pathwise_error} in terms of the continuous loss in Equation~\eqref{eq:mild_loss}.
\begin{theorem}[Continuous-loss control of the rollout]\label{thm:theory_loss}
Work in the setting of Section~\ref{sec:theory_setting}, with deterministic $L_a,L_z$, exact initial data and paths, and measurable $C^1([0,H];\mathbb R^K)$-valued residuals with finite $J,\overline{\mathcal D}_R$. Then
\begin{equation}
\label{eq:theory_loss_error}
 \left(\mathbb E[\|a-\widehat a\|_{\infty,[0,T]}^2]\right)^{1/2} \le 2Te^{L_aT}\left(J^{1/2}+\overline{\mathcal D}_R^{1/3}J^{1/3}\right).
\end{equation}
No independence between residuals or their regularity bounds is required.
\end{theorem}
This result provides a conditional accuracy guarantee for BRiSE's continuous-loss formulation. The proof is given in Section~\ref{sec:theory_decoder}; extending the guarantee to empirical training requires quadrature and sampling control.

\section{Further Theoretical Analysis} \label{sec:theory}

This section proves the decoder, residual, and loss guarantees in Section~\ref{sec:method}, then extends them to discrete training and temporal--spatial errors. All paths share the same Wiener realization. The analysis uses contractive semigroup estimates \cite{hochbruck2010exponential,jentzen2011efficient} and stability of the BRiSE mild residual.

\subsection{Analytical Setting}\label{sec:theory_setting}
The error estimates compare two trajectories in the same Galerkin space. We first specify the paths and the region on which their drifts can be compared, then record the endpoint update that links successive intervals.

Fix a spectral dimension $K$, a final time $T>0$, and a uniform partition with $H=T/N$. Let $P_K$ denote the orthogonal projection onto $V_K=\operatorname{span}\{\phi_1,\ldots,\phi_K\}$, and let $I_Kx=\sum_{k=1}^Kx_k\phi_k$. By orthonormality of the Galerkin basis, $\|I_Kx\|=\|x\|_2$.
The exact Galerkin remainder $a$ satisfies $\dot a=-\Lambda a+G(a,z)$, with $a(0)=I_K^{-1}P_Ku_0$. The predicted interval states are denoted by $x_n$. We work with $\Lambda=\operatorname{diag}(\lambda_1,\ldots,\lambda_K)$ with $0\le\lambda_1\le\cdots\le\lambda_K$. The paths $z,\widetilde z$ belong to $C([0,T];\mathbb R^K)$, and the exact Galerkin solution exists on $[0,T]$. Here $\widetilde z$ is the path used for numerical reconstruction and residual evaluation; $\widetilde z=z$ corresponds to exact path information.

On the states compared below, the nonlinear drift admits finite constants $L_a,L_z\ge0$ such that
\begin{equation}
\label{eq:theory_lipschitz} \|G(x,y)-G(\bar x,\bar y)\|_2 \le L_a\|x-\bar x\|_2+L_z\|y-\bar y\|_2
\end{equation}
for all pairs of states compared in the analysis. The bound in Equation~\eqref{eq:theory_lipschitz} may hold globally or on a specified bounded region containing the relevant trajectories. In the latter case, the conclusions apply only while those trajectories remain in that region. For the semilinear heat benchmark, with $A=\partial_{xx}$ and $F(u)=-u-\sin u$, one can take $L_a=L_z=2$ independently of $K$. Polynomial reaction terms and derivative nonlinearities generally require localization or a different function-space analysis.

For Allen--Cahn, Fisher--KPP, and Gray--Scott, we use local Lipschitz bounds on the compared trajectories. At fixed $K$, spectral discretization also gives local coefficient-space bounds for the derivative nonlinearities in Burgers and Cahn--Hilliard; their constants may depend on $K$. For a trained network, population expectations are taken over independent evaluation paths, with network parameters fixed or conditional on the realized training data.

The decoder and residual use the same interval coefficients and realized path. Endpoints are propagated exactly by the decoder, with no additional interface jump. Writing Equation~\eqref{eq:brise_decoder} in diagonal matrix form and evaluating it at $s=H$ gives
\begin{equation}
\label{eq:theory_decoder} \begin{aligned}
v_n(s)&=E_sx_n+D_sG(x_n,\widetilde z(t_n)) +\sum_{j=1}^M B_j(s)c_{n,j},\\
x_{n+1}&=v_n(H),\quad B_j(s)=\operatorname{diag}(\Psi_{j,1}(s;H),\ldots,\Psi_{j,K}(s;H)).
\end{aligned}
\end{equation}
Here $c_{n,j}=c_{\theta,n,j}$ and $M=2$ in BRiSE. The global remainder approximation is $\widehat a(t_n+s)=v_n(s)$. The analysis is pathwise and allows coefficients to depend on within-interval innovations.

\subsection{Decoder Stability and Residual Error Bounds}\label{sec:theory_decoder}
We begin with a single decoder interval to identify how input perturbations affect the reconstructed trajectory. The two rollout estimates of Section~\ref{sec:pathwise_learning} then require two further steps. We first track how the interval residuals accumulate along a path, then bound their accumulated size by the continuous training loss.

In the setting of Section~\ref{sec:theory_setting}, the decoder in Equation~\eqref{eq:theory_decoder} has the following consistency and stability properties. Here $v_n=\widehat a_{\theta,n}$ in the notation of Section~\ref{sec:method}.

\begin{lemma}[Decoder consistency and perturbation stability]\label{lem:theory_decoder}
The decoder satisfies $v_n(0)=x_n$ and $v_n'(0)=-\Lambda x_n+G(x_n,\widetilde z(t_n))$. For two interval decoders with the same $H,\Lambda$ and data $(x,y,c_j)$ and $(\bar x,\bar y,\bar c_j)$, where $y$ denotes the initial noise state and $c_j$ the forcing coefficients,
\begin{equation}
\label{eq:theory_decoder_stability} \begin{aligned}
\|v-\bar v\|_{\infty,[0,H]} &\le (1+HL_a)\|x-\bar x\|_2+ HL_z\|y-\bar y\|_2 + H\sum_{j=1}^M\frac{\|c_j-\bar c_j\|_2}{j+1}.
\end{aligned}
\end{equation}
The numerical basis constants are independent of $K$ and $\lambda_K$.
\end{lemma}

The decoder is continuously differentiable. Using the operator norm induced by the Euclidean vector norm, we obtain the following bounds for the semigroup and basis matrices in Lemma~\ref{lem:theory_decoder}. We have
\begin{equation}
\label{eq:theory_basis_bounds} \|E_s\|_{\mathrm{op}}\le1,\quad \|D_s\|_{\mathrm{op}}\le s,\quad \|B_j(s)\|_{\mathrm{op}}\le\frac{s^{j+1}}{(j+1)H^j}, \quad 0\le s\le H.
\end{equation}
These bounds include zero eigenvalues.
\begin{proof}[Proof of Lemma~\ref{lem:theory_decoder}]\label{proof:theory_decoder}
Every diagonal entry of $E_s$ lies in $[0,1]$. Consequently,
$$ \|D_s\|_{\mathrm{op}}\le\int_0^s1\,\mathrm{d}r=s,\quad \|B_j(s)\|_{\mathrm{op}}\le\int_0^s(r/H)^j\,\mathrm{d}r =\frac{s^{j+1}}{(j+1)H^j}. $$
These bounds also hold for zero eigenvalues. Differentiating the integral representations of $D_s$ and $B_j(s)$ with respect to the local time $s$, with the step size $H$ held fixed, gives
$$ D_s'=-\Lambda D_s+I,\quad B_j'(s)=-\Lambda B_j(s)+(s/H)^jI. $$
Substituting these derivatives into the derivative of Equation~\eqref{eq:theory_decoder} gives
\begin{equation}
\label{eq:theory_forcing_identity} v_n'(s)=-\Lambda v_n(s)+f_n(s),\quad f_n(s)=G(x_n,\widetilde z(t_n))+\sum_{j=1}^M(s/H)^jc_{n,j}.
\end{equation}
At $s=0$, the integral and nonconstant forcing terms vanish, proving consistency. Subtraction and Equation~\eqref{eq:theory_lipschitz} give
$$ \|v(s)-\bar v(s)\|_2 \le(1+sL_a)\|x-\bar x\|_2+sL_z\|y-\bar y\|_2 +\sum_{j=1}^M\frac{s^{j+1}}{(j+1)H^j}\|c_j-\bar c_j\|_2. $$
Taking the supremum over $0\le s\le H$ proves Equation~\eqref{eq:theory_decoder_stability}.
\end{proof}

If the stacked coefficient network is $L_N$-Lipschitz in the Euclidean norms of the unscaled input and stacked output, Cauchy--Schwarz gives
\begin{equation}
\label{eq:theory_network_sensitivity} \sum_{j=1}^M\frac{\|c_j-\bar c_j\|_2}{j+1} \le\beta_M L_N\|\boldsymbol{\chi}-\bar{\boldsymbol{\chi}}\|_2,\quad \beta_M=\left(\sum_{j=1}^M(j+1)^{-2}\right)^{1/2}.
\end{equation}
For identical noise inputs, a perturbation of the starting remainder state thus propagates by at most $1+H(L_a+\beta_ML_N)$ per step, or $\exp((L_a+\beta_ML_N)T)$ over the horizon. This sensitivity constant differs from the residual stability constant in Theorem~\ref{thm:theory_residual}. For a feedforward network, a sufficient Lipschitz bound is the product of the weight operator norms and activation Lipschitz constants, multiplied by the norms of any input and output scaling operators.

For two OU paths with supremum distance $\varepsilon_z$ on an interval, their innovations differ by at most $2\varepsilon_z$, their endpoint innovations by at most $2\varepsilon_z$, and their exact integrated innovations by at most $2H\varepsilon_z$. An error in the numerical integrated feature contributes an additional input perturbation through Equation~\eqref{eq:theory_network_sensitivity}.

\paragraph{Residual error and mean-square extension.}\phantomsection\label{sec:theory_residual}
We now follow the residuals from one interval to the next to bound the error along a single path, then take a mean-square norm over paths. Define $e_0=a(0)-x_0$, $\rho_n=\|R_n\|_{\infty,[0,H]}$, and the continuous mild residual $R_n=\mathcal R_{\theta,n}$ by
\begin{equation}
\label{eq:theory_residual} R_n(s)={}v_n(s)-E_sx_n-\int_0^s E_{s-r}G(v_n(r),\widetilde z(t_n+r))\,\mathrm{d}r.
\end{equation}
We now have the ingredients for the proof of Theorem~\ref{thm:theory_residual}. Contractivity of $E_s$ and the drift bound control the error within one interval. The identity $x_{n+1}=v_n(H)$ carries its endpoint error into the next interval.
\begin{proof}[Proof of Theorem~\ref{thm:theory_residual}]\label{proof:theory_residual}
Solve Equation~\eqref{eq:theory_residual} for $v_n(s)$ and subtract it from the exact mild solution in Equation~\eqref{eq:interval_mild}. Let
$$ \varepsilon_n(s)=\|a(t_n+s)-v_n(s)\|_2,\quad d_n=\int_0^H\|z(t_n+r)-\widetilde z(t_n+r)\|_2\,\mathrm{d}r. $$
Contractivity and Equation~\eqref{eq:theory_lipschitz} imply
\begin{equation}
\label{eq:theory_local_volterra} \varepsilon_n(s)\le \varepsilon_n(0)+\rho_n+L_zd_n+L_a\int_0^s\varepsilon_n(r)\,\mathrm{d}r.
\end{equation}
Writing $b_n=\varepsilon_n(0)+\rho_n+L_zd_n$, iteration of this scalar integral inequality gives
$$ \varepsilon_n(s)\le b_n\sum_{\ell=0}^{m-1}\frac{(L_as)^\ell}{\ell!} +\left(\sup_{r\in[0,H]}\varepsilon_n(r)\right)\frac{(L_as)^m}{m!}. $$
The final term tends to zero as $m\to\infty$, since $\varepsilon_n$ is continuous on a compact interval. Hence
$$ \sup_{r\le s}\varepsilon_n(r)\le e^{L_as}(\varepsilon_n(0)+\rho_n+L_zd_n). $$
By the endpoint update in Equation~\eqref{eq:theory_decoder}, the error at the end of one interval is the initial error on the next, $\varepsilon_n(H)=\varepsilon_{n+1}(0)$. Iterating the interval estimate at these shared endpoints gives
$$ \varepsilon_n(0)\le e^{L_at_n}\|e_0\|_2 +\sum_{i<n}e^{L_a(t_n-t_i)}(\rho_i+L_zd_i). $$
Combining the last two inequalities and bounding each propagation factor by $e^{L_aT}$ proves Equation~\eqref{eq:theory_pathwise_error}.
\end{proof}

Across $N=T/H$ intervals, local defects of size $H\varepsilon$ thus contribute at most $T\varepsilon$ before multiplication by the stability factor $e^{L_aT}$. This estimate requires no additional restriction on $H\lambda_K$.

If $L_a,L_z$ are deterministic and the terms on the right have finite second moments, taking the $L^2(\Omega)$ norm in Equation~\eqref{eq:theory_pathwise_error} and applying Minkowski's inequality to the sum and the integral gives
\begin{equation}
\label{eq:theory_rms_error} \big(\mathbb E[\|a-\widehat a\|_{\infty,[0,T]}^2]\big)^{1/2} \le e^{L_aT}\bigg(\|e_0\|_{L^2(\Omega)}+\sum_{n=0}^{N-1}\|\rho_n\|_{L^2(\Omega)}+L_z\int_0^T\|z(t)-\widetilde z(t)\|_{L^2(\Omega)}\,\mathrm{d}t\bigg).
\end{equation}

For locally Lipschitz nonlinearities, the proof applies while all compared trajectories remain in the region on which Equation~\eqref{eq:theory_lipschitz} holds.

\paragraph{Continuous-loss control.}\phantomsection\label{sec:theory_loss}
We now use the population training loss to bound the residual contribution in Theorem~\ref{thm:theory_residual} \cite{mishra2023estimates}. Define
$$
J=\mathbb E\!\Bigg[\frac1{NH^3}\sum_{n=0}^{N-1} \int_0^H\|R_n(s)\|_2^2\,\mathrm{d}s\Bigg], \quad
\overline{\mathcal D}_R =\left(\mathbb E\Bigg[\max_n\sup_{s\in[0,H]}\|R_n'(s)\|_2^2\Bigg]\right)^{1/2}.
$$
Here $J$ denotes the population loss, distinct from the integrated innovation $J_n$ in Section~\ref{subsec:OU Bridge-Informed Path Encoding}; for exact paths, $J=\mathcal L_{\mathrm{mild}}(\theta)$. Expectations are over evaluation paths for fixed trained parameters. For one noise realization $\omega$, define its normalized loss by
\begin{equation}
\label{eq:theory_population_loss} J_\omega=\frac1N\sum_{n=0}^{N-1}\frac1H\int_0^H \frac{\|R_n(s)\|_2^2}{H^2}\,\mathrm{d}s,\quad J=\mathbb E[J_\omega].
\end{equation}
The residual bound underlying Theorem~\ref{thm:theory_loss} also holds for inexact initial data and paths. Work in the setting of Section~\ref{sec:theory_setting}, with deterministic $L_a,L_z$. Suppose the residuals are measurable $C^1([0,H];\mathbb R^K)$-valued random variables and
$$ \mathcal D_R=\max_{0\le n<N}\sup_{s\in[0,H]}\|R_n'(s)\|_2, \quad \overline{\mathcal D}_R=\|\mathcal D_R\|_{L^2(\Omega)}<\infty, \quad J<\infty. $$
Then
\begin{equation}
\label{eq:theory_loss_rho} \left\|\sum_{n=0}^{N-1}\rho_n\right\|_{L^2(\Omega)} \le 2T\left(J^{1/2}+\overline{\mathcal D}_R^{\,1/3}J^{1/3}\right).
\end{equation}
Equation~\eqref{eq:theory_loss_rho} controls the sum of the interval suprema and allows $\mathcal D_R$ and $J_\omega$ to be dependent. For exact initial data and paths, this is the only error term in Equation~\eqref{eq:theory_pathwise_error}, giving Equation~\eqref{eq:theory_loss_error}. With bounded $L_a$, the conditions $J\to0$ and $\overline{\mathcal D}_R J\to0$ yield convergence through Equation~\eqref{eq:theory_loss_error}. Section~\ref{sec:theory_discrete} translates these conditions into empirical-residual criteria.

A small time-averaged loss can hide a short residual peak. A bound on the residual derivative limits how narrow a peak of a given height can be. The following lemma uses this observation to bound the peak height from the time average of the squared residual.

\begin{lemma}[A temporal interpolation estimate]\label{lem:theory_interpolation}
Let $f:[0,1]\to\mathbb R^K$ be $L_f$-Lipschitz, and let $l=\|f\|_{L^2(0,1;\mathbb R^K)}$. Then
\begin{equation}
\label{eq:theory_interpolation} \|f\|_{\infty,[0,1]}\le2l+2L_f^{1/3}l^{2/3}.
\end{equation}
The numerical constants are independent of $K$.
\end{lemma}
\begin{proof}
If $L_f=0$, the function is constant. Otherwise write $M_f=\|f\|_{\infty,[0,1]}$. The case $M_f=0$ is immediate. At a maximizer $\tau_*$, the reverse triangle inequality implies $\|f(\tau)\|_2\ge M_f/2$ whenever $|\tau-\tau_*|\le M_f/(2L_f)$. Intersecting this neighborhood with $[0,1]$ gives a set of length at least $\min\{1,M_f/(2L_f)\}$. Consequently,
$$ l^2\ge\frac{M_f^2}{4}\min\left\{1,\frac{M_f}{2L_f}\right\}. $$
If the minimum is 1, then $M_f\le2l$. Otherwise $M_f^3\le8L_fl^2$. Taking the sum of the two resulting bounds proves Equation~\eqref{eq:theory_interpolation}.
\end{proof}

Lemma~\ref{lem:theory_interpolation} completes the preparation for the proof of Theorem~\ref{thm:theory_loss}. After rescaling each interval to $[0,1]$, we can apply the lemma to its residual and sum the resulting bounds in the pathwise estimate of Theorem~\ref{thm:theory_residual}.

\begin{proof}[Proof of Theorem~\ref{thm:theory_loss}]\label{proof:theory_loss}
Set $\mathfrak r_n(\tau)=R_n(H\tau)/H$ and $l_n^2=\int_0^1\|\mathfrak r_n(\tau)\|_2^2\,\mathrm{d}\tau$. The Lipschitz constant of $\mathfrak r_n$ is bounded by $\mathcal D_R$, since $\mathfrak r_n'(\tau)=R_n'(H\tau)$. Lemma~\ref{lem:theory_interpolation} implies
$$ \rho_n\le2H(l_n+\mathcal D_R^{1/3}l_n^{2/3}). $$
The change of variables $s=H\tau$ in Equation~\eqref{eq:theory_population_loss} gives $J_\omega=N^{-1}\sum_n l_n^2$. Cauchy--Schwarz and concavity then give
$$ \frac1N\sum_n l_n\le J_\omega^{1/2},\quad \frac1N\sum_n l_n^{2/3}\le J_\omega^{1/3}. $$
Using $NH=T$ yields the pathwise inequality
\begin{equation}
\label{eq:theory_pathwise_loss} \sum_n\rho_n\le2T(J_\omega^{1/2}+\mathcal D_R^{1/3}J_\omega^{1/3}).
\end{equation}
For the mixed term, H\"older's inequality with exponents $3$ and $3/2$ gives
$$ \mathbb E[\mathcal D_R^{2/3}J_\omega^{2/3}] \le(\mathbb E[\mathcal D_R^2])^{1/3}(\mathbb E[J_\omega])^{2/3}. $$
Take square roots and apply the $L^2(\Omega)$ triangle inequality to Equation~\eqref{eq:theory_pathwise_loss} to obtain Equation~\eqref{eq:theory_loss_rho}. For Equation~\eqref{eq:theory_loss_error}, take the $L^2(\Omega)$ norm of Equation~\eqref{eq:theory_pathwise_error}, keeping $\sum_n\rho_n$ together as a single random variable, and then apply Equation~\eqref{eq:theory_loss_rho}.
\end{proof}

Together, Theorems~\ref{thm:theory_residual} and~\ref{thm:theory_loss} justify the choice of the mild residual as the training objective in Section~\ref{sec:pathwise_learning}. The first theorem shows how interval residuals propagate into the trajectory error. The second shows that, for exact initial data and noise paths and a uniform bound on $\overline{\mathcal D}_R$, a small population residual loss guarantees a small root-mean-square rollout error for the Galerkin solution.

The residual regularity assumption can be expressed in terms of the actual decoder. Inserting the mild solution of Equation~\eqref{eq:theory_forcing_identity} for $v_n(s)$ in Equation~\eqref{eq:theory_residual} cancels the common term $E_sx_n$ and gives
$$ R_n(s)=\int_0^sE_{s-r}\delta_n(r)\,\mathrm{d}r,\quad \delta_n(r)=f_n(r)-G(v_n(r),\widetilde z(t_n+r)). $$
Since $\widetilde z$ is continuous,
$$ R_n'=-\Lambda R_n+\delta_n,\quad \mathcal D_R\le(1+\lambda_KH)\max_n\|\delta_n\|_{\infty,[0,H]}. $$
If $\|G(x_n,\widetilde z(t_n))\|_2\le B_0$, $\|c_{n,j}\|_2\le C_j$, and $\|G(v_n(s),\widetilde z(t_n+s))\|_2\le B_F$ uniformly in $n,s$, with square-integrable random bounds, then
\begin{equation}
\label{eq:theory_regularity_constant} \overline{\mathcal D}_R \le(1+\lambda_KH) \left(\|B_0\|_{L^2(\Omega)}+\sum_j\|C_j\|_{L^2(\Omega)} +\|B_F\|_{L^2(\Omega)}\right).
\end{equation}
Equation~\eqref{eq:theory_regularity_constant} controls residual regularity through the forcing defect and $\lambda_KH$, without differentiating the OU path. Uniform refinement estimates require the displayed bound to remain controlled.

\subsection{Discrete Training and Convergence Estimates}\label{sec:theory_discrete}
The continuous-loss estimate involves an expectation and a time integral, whereas training uses finitely many paths and collocation times. We first quantify this discrepancy and derive a convergence criterion from the sampled loss. We then obtain a temporal estimate at fixed $K$ from OU increments and the sizes of the learned coefficients. Finally, we include spatial truncation and identify which constants must be controlled when the spatial and temporal resolutions vary together.

For a single realized path, let $\widehat R_n(s_m)$ denote the quadrature residual $\widehat{\mathcal R}_{\theta,n}(s_m)$ in Appendix~\ref{app:numerical_residual_evaluation}, using the path $\widetilde z$. Define
$$ \delta_c=\sup_{s\in[0,H]}\min_m|s-s_m|,\quad \varepsilon_{\mathrm{quad},n}=\max_m\|R_n(s_m)-\widehat R_n(s_m)\|_2,\quad \ell_n=\frac{1}{N_{\mathrm c}}\sum_{m=1}^{N_{\mathrm c}}\frac{\|\widehat R_n(s_m)\|_2^2}{H^2}. $$
Comparison with the nearest collocation point gives
\begin{equation}
\label{eq:theory_discrete_certificate} \begin{aligned}
&\rho_n\le H\sqrt{N_{\mathrm c}\ell_n}+\varepsilon_{\mathrm{quad},n}+\mathcal D_R\delta_c,\\
&\sum_n\rho_n\le T\sqrt{N_{\mathrm c}\widehat J_\omega} +\sum_n\varepsilon_{\mathrm{quad},n}+N\mathcal D_R\delta_c,\quad \widehat J_\omega=\frac1N\sum_n\ell_n.
\end{aligned}
\end{equation}
Indeed, the largest sampled residual is at most the square root of the sum of its squared values, and Cauchy--Schwarz bounds $N^{-1}\sum_n\sqrt{\ell_n}$ by $\sqrt{\widehat J_\omega}$.

The term $\mathcal D_R\delta_c$ bounds how much a residual peak can exceed the values at the collocation times. The quadrature error $\varepsilon_{\mathrm{quad},n}$ arises when evaluating the residual at those times. To see how this error depends on the quadrature mesh, consider composite left quadrature with maximum panel size $\delta_q$. Suppose $g_n(r)=G(v_n(r),\widetilde z(t_n+r))$ satisfies $\|g_n\|_{\infty,[0,H]}\le B_g$ and $[g_n]_{C^\alpha([0,H])}\le C_g$, where $[g]_{C^\alpha(I)}=\sup_{r,r'\in I,\,r\ne r'}\|g(r)-g(r')\|_2/|r-r'|^\alpha$ is the H\"older seminorm for $0<\alpha\le1$. For $r,r'\in[0,s]$,
$$ \|E_{s-r}g_n(r)-E_{s-r'}g_n(r')\|_2 \le C_g|r-r'|^\alpha+B_g\lambda_K|r-r'|. $$
Integrating this difference on each quadrature panel proves
\begin{equation}
\label{eq:theory_quadrature} \varepsilon_{\mathrm{quad},n}\le H(C_g\delta_q^\alpha+B_g\lambda_K\delta_q).
\end{equation}
Comparison with both panel endpoints gives the same bound for composite trapezoidal quadrature. For the exponential product integration rule in Appendix~\ref{app:numerical_residual_evaluation}, the semigroup kernel is integrated exactly against the piecewise-linear interpolant of $g_n$. Contractivity and the interpolation error bound $C_g\delta_q^\alpha$ then give $\varepsilon_{\mathrm{quad},n}\le H C_g\delta_q^\alpha$.

The population loss $J(\theta)$ uses the exact mild integral in the residual. Replacing that integral by numerical quadrature gives $J_Q(\theta)$, which still averages over the path and time distributions. Averaging the same numerical residual loss over the sampled paths and collocation times gives the empirical loss $\widehat J(\theta)$. If the normalized residual quadrature error has population $L^2$ norm at most $\varepsilon_q$, the triangle inequality gives
$$ |\sqrt{J(\theta)}-\sqrt{J_Q(\theta)}|\le\varepsilon_q. $$
Let $\Theta$ be the admissible parameter set of the shared coefficient network and $\varepsilon_{\rm opt}\ge0$ the empirical optimization tolerance. Assume that $|J_Q(\theta)-\widehat J(\theta)|\le g$ for all $\theta\in\Theta$ and that the returned parameter $\widehat\theta\in\Theta$ satisfies
$$ \widehat J(\widehat\theta)\le\inf_{\theta\in\Theta}\widehat J(\theta) +\varepsilon_{\rm opt}. $$
For any comparator $\theta^*\in\Theta$,
$$ J_Q(\widehat\theta) \le\widehat J(\widehat\theta)+g \le\widehat J(\theta^*)+\varepsilon_{\rm opt}+g \le J_Q(\theta^*)+2g+\varepsilon_{\rm opt}. $$
Taking square roots and applying the quadrature bound twice yields
\begin{equation}
\label{eq:theory_risk_decomposition} \sqrt{J(\widehat\theta)} \le\sqrt{J(\theta^*)}+2\varepsilon_q+\sqrt{2g+\varepsilon_{\rm opt}}.
\end{equation}
The comparator uses the shared coefficient network. Its error includes temporal forcing approximation, OU encoding, and feature-to-coefficient approximation; the estimate assumes the stated optimization and uniform risk gaps.

Equation~\eqref{eq:theory_risk_decomposition} bounds the population loss at a given step size. Applying Theorem~\ref{thm:theory_loss} as $H$ decreases turns that bound into a criterion for temporal convergence.

\begin{corollary}[A residual-loss criterion for temporal convergence]\label{cor:theory_training_rate}
Fix $K,T$ and $0<H\le\min\{1,T\}$ with $T/H\in\mathbb N$. Assume the hypotheses of Theorem~\ref{thm:theory_loss} for each $H$, with exact initial data and paths and $L_a\le L_*$ independently of $H$. Fix the trained parameters, or condition on training data. Suppose
$$ \sqrt{J_H}\le\eta_H,\quad \overline{\mathcal D}_{R,H}\le C_DH^\beta,\quad \eta_H\le C_\eta H^q, $$
where $C_D,C_\eta$ are finite deterministic constants in this conditional setting, independent of $H$, $q>0$, and $\beta+2q>0$. Then
\begin{equation}
\label{eq:theory_training_rate} \big(\mathbb E[\|a-\widehat a_H\|_{\infty,[0,T]}^2]\big)^{1/2} \le2Te^{L_*T}\left(\eta_H+C_D^{1/3}H^{\beta/3}\eta_H^{2/3}\right) \le2Te^{L_*T}\left(C_\eta H^q+ C_D^{1/3}C_\eta^{2/3}H^{(\beta+2q)/3}\right).
\end{equation}
Thus the guaranteed rate is $p=\min\{q,(\beta+2q)/3\}$. In particular, a uniformly bounded $\overline{\mathcal D}_{R,H}$ and $J_H=O(H^{3/2})$ suffice for a root-mean-square error of order $H^{1/2}$.

For an empirical certificate, it suffices that the normalized quadrature error is bounded by $\varepsilon_{q,H}$ and $|J_{Q,H}(\widehat\theta_H)-\widehat J_H(\widehat\theta_H)|\le g_H$ for the returned model. One may then take
\begin{equation}
\label{eq:theory_empirical_certificate} \eta_H=\varepsilon_{q,H}+ \sqrt{\widehat J_H(\widehat\theta_H)+g_H}.
\end{equation}
Alternatively, under the uniform risk and optimization hypotheses leading to Equation~\eqref{eq:theory_risk_decomposition}, one may take $\eta_H=\sqrt{J_H(\theta_H^*)}+2\varepsilon_{q,H}+ \sqrt{2g_H+\varepsilon_{{\rm opt},H}}$.
\end{corollary}
\begin{proof}
Theorem~\ref{thm:theory_loss} bounds the error by $2Te^{L_*T}(\sqrt{J_H}+\overline{\mathcal D}_{R,H}^{1/3}J_H^{1/3})$. By monotonicity in $J_H$, substitute $\sqrt{J_H}\le\eta_H$ and then the two assumed power bounds to obtain Equation~\eqref{eq:theory_training_rate}. Since $H\le1$, each power is bounded above by $H^p$ with the stated $p$. For Equation~\eqref{eq:theory_empirical_certificate}, the population triangle inequality and the risk-gap bound give
$$ \sqrt{J_H(\widehat\theta_H)} \le\varepsilon_{q,H}+\sqrt{J_{Q,H}(\widehat\theta_H)} \le\varepsilon_{q,H}+\sqrt{\widehat J_H(\widehat\theta_H)+g_H}. $$
The alternative choice follows from Equation~\eqref{eq:theory_risk_decomposition}.
\end{proof}

For example, when $\beta=0$, sufficient empirical conditions for the half-order bound are $\varepsilon_{q,H}=O(H^{3/4})$ and $\widehat J_H(\widehat\theta_H)+g_H=O(H^{3/2})$, together with the uniform residual-derivative and stability bounds. If the risk certificate is conditional on a training-data event, the estimate holds on that event.

\paragraph{Conditional temporal convergence.}\phantomsection\label{sec:theory_temporal}
We next derive a temporal convergence rate at fixed spectral dimension by estimating the local forcing error from the OU increments and the predicted correction coefficients. This complements the rate in Corollary~\ref{cor:theory_training_rate}, which is expressed in terms of residual decay.

At fixed $K$, let $Q_K=\operatorname{diag}(q_1,\ldots,q_K)$ be the covariance restricted to the retained modes, expressed in coefficient coordinates. The OU increments satisfy
\begin{equation}
\label{eq:theory_ou_increment} \mathbb E[\|z(t+s)-z(t)\|_2^2] \le 2\sigma^2\operatorname{Tr}(Q_K)s,\quad \operatorname{Tr}(Q_K)=\sum_{k=1}^Kq_k.
\end{equation}
For moment regularity, see \cite{kruse2012regularity}. For each interval, let $a_n^*$ be the exact local remainder solution driven by $z$ and initialized at the predicted state $x_n$. The next theorem combines OU increment regularity with a correction-coefficient bound to give BRiSE a strong temporal convergence estimate.

\begin{theorem}[Conditional strong temporal estimate]\label{thm:theory_temporal}
Fix $K,T$ and let $H=T/N\downarrow0$. Let $z$ be the zero-initial OU process defined in Section~\ref{sec:preliminaries}, with independent Brownian coordinates, finite $\operatorname{Tr}(Q_K)$, and deterministic noise amplitude $\sigma$. Work in the setting of Section~\ref{sec:theory_setting} with deterministic, $H$-independent $L_a,L_z$, exact path information $\widetilde z=z$, and $x_0=a(0)$. Suppose the local solutions $a_n^*\in C^1([0,H];\mathbb R^K)$ exist and the measurable coefficient vectors satisfy
\begin{equation}
\label{eq:theory_temporal_conditions} \sup_{n,s}\|(a_n^*)'(s)\|_{L^2(\Omega)}\le B_a,\quad \sum_{j=1}^M\frac{\|c_{n,j}\|_{L^2(\Omega)}}{j+1}\le C_cH^{1/2},
\end{equation}
uniformly in $n,H$, with finite deterministic $B_a,C_c$. Then
\begin{equation}
\label{eq:theory_temporal_bound} \begin{aligned}
\big(\mathbb E[\|a-\widehat a\|_{\infty,[0,T]}^2]\big)^{1/2} \le Te^{L_aT}(1+L_aH)\bigg( \frac{L_aB_a}{2}H +\bigg(\frac{2\sqrt2}{3}L_z|\sigma| \sqrt{\operatorname{Tr}(Q_K)}+C_c\bigg)H^{1/2}\bigg).
\end{aligned}
\end{equation}
For fixed $K$, the root-mean-square temporal error is $O(H^{1/2})$, with the dependence on $T,L_a,L_z,B_a,C_c,\sigma,Q_K$ given by the displayed bound.
\end{theorem}

Theorem~\ref{thm:theory_temporal} gives temporal convergence at fixed spectral resolution under the coefficient bound in Equation~\eqref{eq:theory_temporal_conditions}. Its proof follows below. The zero-correction choice $c_{n,j}=0$ also satisfies Equation~\eqref{eq:theory_temporal_conditions}. Corollary~\ref{cor:theory_training_rate} provides a complementary error estimate directly from the population or certified empirical residual.

We first derive Equation~\eqref{eq:theory_ou_increment}. The exact OU transition decomposes an increment into $(E_s-I)z(t)$ and an independent innovation. Thus
$$ \mathbb E[\|z(t+s)-z(t)\|_2^2] =\sigma^2\sum_{k=1}^Kq_k \bigg((1-e^{-\lambda_ks})^2b_k(t)^2+b_k(s)^2\bigg), $$
where $b_k$ is defined in Appendix~\ref{app:OU Path Construction}. For $\lambda_k>0$, use $b_k(t)^2\le(2\lambda_k)^{-1}$, $(1-e^{-x})^2\le2x$, and $b_k(s)^2\le s$. A zero mode contributes exactly $\sigma^2q_ks$. Summation proves Equation~\eqref{eq:theory_ou_increment}.

\begin{proof}[Proof of Theorem~\ref{thm:theory_temporal}]
Let $g_n^*(s)=G(a_n^*(s),z(t_n+s))$. Since $a_n^*(0)=x_n$, its initial forcing agrees with that in Equation~\eqref{eq:theory_forcing_identity}. Minkowski's inequality, Equation~\eqref{eq:theory_temporal_conditions}, and Equation~\eqref{eq:theory_ou_increment} imply
$$ \|g_n^*(s)-g_n^*(0)\|_{L^2(\Omega)} \le L_aB_as+\sqrt2 L_z|\sigma|\sqrt{\operatorname{Tr}(Q_K)}s^{1/2}. $$
The exact local mild equation and the decoder formula give, pathwise,
$$ \|a_n^*-v_n\|_{\infty,[0,H]} \le\int_0^H\|g_n^*(r)-g_n^*(0)\|_2\,\mathrm{d}r +H\sum_{j=1}^M\frac{\|c_{n,j}\|_2}{j+1}. $$
Taking the $L^2(\Omega)$ norm and integrating the preceding increment estimate yields
$$ \big(\mathbb E[\|a_n^*-v_n\|_{\infty,[0,H]}^2]\big)^{1/2} \le \frac{L_aB_a}{2}H^2+ \left(\frac{2\sqrt2}{3}L_z|\sigma|\sqrt{\operatorname{Tr}(Q_K)} +C_c\right)H^{3/2}. $$
To connect this local comparison with the residual, subtract the exact local mild formula from Equation~\eqref{eq:theory_residual}:
$$ R_n(s)=v_n(s)-a_n^*(s) +\int_0^sE_{s-r} \left[G(a_n^*(r),z(t_n+r))-G(v_n(r),z(t_n+r))\right]\,\mathrm{d}r. $$
Consequently,
$$ \|\rho_n\|_{L^2(\Omega)} \le(1+L_aH) \big(\mathbb E[\|a_n^*-v_n\|_{\infty,[0,H]}^2]\big)^{1/2}. $$
Sum over $N=T/H$ intervals and apply Equation~\eqref{eq:theory_rms_error} with zero initial and path errors. This gives Equation~\eqref{eq:theory_temporal_bound}.
\end{proof}

More generally, Equation~\eqref{eq:theory_pathwise_error} implies an order-$p$ pathwise bound if $\rho_n\le C_rH^{p+1}$, $\|e_0\|_2\le C_0H^p$, and $\|z-\widetilde z\|_{\infty,[0,T]}\le C_zH^p$, with constants independent of $H$:
$$ \|a-\widehat a\|_{\infty,[0,T]} \le e^{L_aT}(C_0+TC_r+L_zTC_z)H^p. $$
For a local forcing polynomial $f_n$, set $\varepsilon_{f,n}=\|g_n^*-f_n\|_{\infty,[0,H]}$. Direct subtraction gives
$$ \|a_n^*-v_n\|_{\infty,[0,H]}\le H\varepsilon_{f,n}, \quad \rho_n\le H(1+L_aH)\varepsilon_{f,n}. $$
For $\alpha$-H\"older forcing, freezing $g_n^*$ gives $\varepsilon_{f,n}\le[g_n^*]_{C^\alpha}H^\alpha$. Finite-dimensional OU paths admit every exponent $\alpha<1/2$. An $O(H^3)$ forcing approximation instead requires three bounded time derivatives and matching coefficient accuracy.

\paragraph{Spatial truncation and total error.}\phantomsection\label{sec:theory_total}
We now pass from the Galerkin trajectory to the full SPDE solution by adding the spatial truncation error to the temporal and noise-path errors in the retained modes. The full reconstruction is $\widehat u=I_K(\widehat a+\widetilde z)$ and the exact Galerkin solution is $u_K=I_Ka+Z_K$, where $Z_K=P_KZ=I_Kz$. Hence
\begin{equation}
\label{eq:theory_total_split} u-\widehat u=(u-u_K)+I_K(a-\widehat a)+I_K(z-\widetilde z).
\end{equation}
Generally $u_K\ne P_Ku$. For a globally Lipschitz map $F:\mathcal H\to\mathcal H$ with constant $L$, define
$$ U_{0,K}=\|(I-P_K)u_0\|,\quad F_K^{\rm tail}=\sup_{t\le T}\|(I-P_K)F(u(t))\|,\quad Z_K^{\rm tail}=\sup_{t\le T}\|(I-P_K)Z(t)\|. $$
For existing continuous solutions, subtracting the continuous and projected remainder equations and using
$$ F(u)-P_KF(u_K)=(I-P_K)F(u)+P_K(F(u)-F(u_K)) $$
gives
$$ \|v(t)-v_K(t)\| \le U_{0,K}+TF_K^{\rm tail}+LTZ_K^{\rm tail} +L\int_0^t\|v(r)-v_K(r)\|\,\mathrm{d}r. $$
Applying Gr\"onwall's inequality to this remainder bound and adding the noise tail $Z_K^{\rm tail}$ using $u-u_K=(v-v_K)+(Z-Z_K)$ gives
\begin{equation}
\label{eq:theory_spatial_tail} \|u-u_K\|_{\infty,[0,T];\mathcal H} \le e^{LT}(U_{0,K}+TF_K^{\rm tail}) +(1+LTe^{LT})Z_K^{\rm tail} =:\mathcal E_{\rm sp}(K).
\end{equation}
Related spectral Galerkin estimates appear in \cite{kruse2014optimal,jentzen2011efficient}. Here the spatial tail follows directly from the preceding inequality. For $r\ge0$ and $w\in D((I-A)^{r/2})$, Parseval's identity gives
$$ \|(I-P_K)w\|^2 =\sum_{k>K}|\langle w,\phi_k\rangle|^2 \le(1+\inf_{k>K}\lambda_k)^{-r}\|(I-A)^{r/2}w\|^2. $$
Applying this inequality to the three tails in Equation~\eqref{eq:theory_spatial_tail} requires the corresponding uniform-in-time regularity norms. A root-mean-square supremum estimate additionally requires their second moments. In contrast, the exact identity
$$ \mathbb E\big[\|(I-P_K)Z(t)\|^2\big] =\sigma^2\sum_{k>K}q_kb_k(t)^2 $$
is only a fixed-time bound and cannot be substituted for the expected temporal supremum. The spatial tail estimate now controls the first term in Equation~\eqref{eq:theory_total_split}. Combining it with the residual estimate for the retained modes gives the following bound for the full reconstruction.

\begin{corollary}[Conditional total error bound]\label{cor:theory_total}
Assume the hypotheses of Theorem~\ref{thm:theory_residual}, the general residual moment and regularity hypotheses of Equation~\eqref{eq:theory_loss_rho}, and the globally Lipschitz spatial setting of Equation~\eqref{eq:theory_spatial_tail}. Suppose $\mathcal E_{\rm sp}(K)$, $e_0$, and $\varepsilon_{\rm path}=\|z-\widetilde z\|_{\infty,[0,T]}$ have finite second moments. Then
\begin{equation}
\label{eq:theory_total_bound} \begin{aligned}
\big(\mathbb E\big[\|u-\widehat u\|_{\infty,[0,T];\mathcal H}^2\big]\big)^{1/2} \le{}&\|\mathcal E_{\rm sp}(K)\|_{L^2(\Omega)} +(1+L_zTe^{L_aT})\|\varepsilon_{\rm path}\|_{L^2(\Omega)}\\
&+e^{L_aT}\left( \|e_0\|_{L^2(\Omega)} +2T\left(J^{1/2}+\overline{\mathcal D}_R^{\,1/3}J^{1/3}\right)\right).
\end{aligned}
\end{equation}
All spectral, network, step-size and noise dependence is retained through the displayed constants and error quantities. If the empirical-risk hypotheses of Section~\ref{sec:theory_discrete} hold, $J^{1/2}$ may be bounded by Equation~\eqref{eq:theory_risk_decomposition}.
\end{corollary}
\begin{proof}
Use the isometry of $I_K$ in Equation~\eqref{eq:theory_total_split}, apply Equation~\eqref{eq:theory_spatial_tail}, and then take $L^2(\Omega)$ norms. Apply the pathwise estimate Equation~\eqref{eq:theory_pathwise_error} while keeping the sum of residuals together, and use Equation~\eqref{eq:theory_loss_rho}. The reconstructed noise contributes $\|\varepsilon_{\rm path}\|_{L^2(\Omega)}$ directly and $L_zTe^{L_aT}\|\varepsilon_{\rm path}\|_{L^2(\Omega)}$ through the remainder equation.
\end{proof}

For a prescribed OU subgrid of width $\delta_z$, linear interpolation of exact nodal values satisfies $\|z-\widetilde z\|_{\infty,[0,T]}\le[z]_{C^\alpha([0,T])}\delta_z^\alpha$ whenever the realized path is $\alpha$-H\"older. The trapezoidal integrated feature satisfies
$$ \|J_n-\widehat J_n\|_2\le H[Y_n]_{C^\alpha([0,H])}\delta_z^\alpha, \quad [Y_n]_{C^\alpha([0,H])}\le[z]_{C^\alpha([0,T])} +\lambda_K\|z(t_n)\|_2H^{1-\alpha}. $$
These follow by integrating the error of the piecewise-linear interpolant. The a posteriori residual accounts for the effect of feature approximation on the predicted coefficients.

\paragraph{Dependence on discretization and network parameters.}
To interpret the total error bound under refinement, we identify which constants depend on the retained modes, the quadrature mesh, and the trained network. The semigroup and basis bounds in Equation~\eqref{eq:theory_basis_bounds} are uniform in the spectrum. The nonlinear constant need not be: for a reaction term $u-u^3$, a bound $\|u\|_{L^\infty(D)}\le M_\infty$ gives $L_a,L_z\le1+3M_\infty^2$. If only a coefficient radius is known, the estimate
$$ \|I_Kx\|_{L^\infty(D)}\le C_K\|x\|_2, \quad C_K=\sup_{\xi\in D}\left(\sum_{k=1}^K|\phi_k(\xi)|^2\right)^{1/2} $$
uses the Euclidean pointwise norm for vector-valued basis functions and the corresponding essential-supremum norm for systems, and introduces spectral dependence. Derivative nonlinearities can introduce additional powers of the largest retained wavenumber. Network weights and normalization affect $L_N$, coefficient moments and $\overline{\mathcal D}_R$. Noise enters through $\sigma^2\operatorname{Tr}(Q_K)$, covariance tails and trajectory moments. Inner quadrature introduces $\lambda_K\delta_q$ in Equation~\eqref{eq:theory_quadrature}. Therefore a joint limit $K\to\infty$, $H\to0$ requires control of all these quantities.

\section{Numerical Experiments} \label{sec:numerical_exp}

\subsection{Experimental Setup and Results}

\paragraph{Baselines, benchmarks and implementation details.}
We compare BRiSE against five baselines, including two physics-informed neural methods: PINN \cite{raissi2019physics} and NeS-PINN \cite{baptista2026solving}; three classical numerical solvers: Euler--Maruyama (EM) \cite{kloeden1977numerical}, Exponential Euler (Exp Euler) \cite{jentzen2009overcoming}, and implicit--explicit Euler (IMEX Euler) \cite{ascher1995implicit}. We consider six SPDE benchmarks described in Section~\ref{sec:SPDE}, which includes stochastic Burgers, stochastic Heat, stochastic Fisher--KPP, stochastic Allen--Cahn, stochastic Cahn--Hilliard, and stochastic Gray--Scott equations. These benchmarks cover linear and nonlinear dynamics, reaction-diffusion systems, and higher-order differential operators. BRiSE and the classical solvers use 20 time intervals across all benchmarks, except for Allen--Cahn, where 100 intervals are employed. All experiments are conducted on a workstation equipped with a single NVIDIA GeForce RTX 4090 GPU (24 GB). More details on the experimental setup are provided in Appendix~\ref{app:reproducibility}.

\paragraph{Evaluation metrics.}
We evaluate the accuracy of each method using the space-time relative $L^2$ error, defined as
\begin{equation}
\mathrm{Rel.}\,L^2 = \frac{\|u_\theta-u_{\mathrm{ref}}\|_{L^2([0,T]\times D)}} {\|u_{\mathrm{ref}}\|_{L^2([0,T]\times D)}}.
\end{equation}
Here $u_\theta$ and $u_{\mathrm{ref}}$ are the predicted and reference solutions for the same Wiener realization. For Gray--Scott, the squared pointwise magnitude in the product-space $L^2$ norm sums the squared errors for both species in the numerator and their squared reference values in the denominator. Both norms are evaluated on a common grid $\{(\tau_i,\xi_j):0\le i\le N_t,\ 1\le j\le N_x\}$, where $N_t$ counts the evaluation time intervals and $N_x$ counts the spatial evaluation points. The times $\tau_i$ specify the evaluation grid, and $t_n$ denotes a coarse integration node. We report the mean and sample standard deviation over eight held-out Wiener paths.

\begin{table*}[t]
\centering
\caption{Space-time relative $L^2$ errors across 6 SPDE benchmarks. Results are reported as the mean ± std over 8 held-out Wiener paths; lower is better. The best mean error in each column is highlighted in \colorbox{bestgreen}{green}, and the second best is \underline{underlined}. Improvement denotes the relative reduction in mean error of BRiSE over the second-best method. Values $>1$, highlighted in \colorbox{failred}{red}, indicate failure to converge.}
\label{tab:main_results}
\resizebox{\textwidth}{!}{
\begin{tabular}{lcccccc}
\toprule
Method
& Burgers
& Heat
& Fisher--KPP
& Allen--Cahn
& Cahn--Hilliard
& Gray--Scott \\
\midrule
PINN
& \underline{8.55e-02 ± 5.94e-03}
& 7.81e-01 ± 5.36e-02
& 2.58e-02 ± 3.31e-03
& 1.56e-01 ± 5.88e-03
& \cellcolor{failred}$>1$
& 6.56e-01 ± 5.37e-03 \\

NeS-PINN
& 1.36e-01 ± 1.25e-02
& 4.95e-01 ± 3.56e-02
& 2.69e-02 ± 1.34e-03
& 9.69e-02 ± 3.91e-03
& 2.53e-01 ± 7.04e-03
& 2.23e-02 ± 1.33e-04 \\

EM
& \cellcolor{failred}$>1$
& 4.72e-02 ± 4.80e-03
& \cellcolor{failred}$>1$
& 4.17e-03 ± 1.15e-04
& \cellcolor{failred}$>1$
& \cellcolor{failred}$>1$ \\

Exp Euler
& 2.20e-01 ± 1.42e-01
& \underline{3.75e-02 ± 3.25e-03}
& \underline{1.02e-02 ± 1.06e-03}
& \underline{4.01e-03 ± 1.37e-04}
& \underline{7.69e-02 ± 1.11e-04}
& \underline{6.06e-03 ± 8.87e-04} \\

IMEX Euler
& 8.97e-02 ± 2.02e-02
& 4.95e-01 ± 3.57e-02
& 1.07e-02 ± 1.04e-03
& 4.16e-03 ± 1.11e-04
& 2.06e-01 ± 6.78e-03
& 1.06e-02 ± 4.69e-04 \\

\midrule
\textbf{BRiSE}
& \cellcolor{bestgreen}7.56e-02 ± 8.52e-03
& \cellcolor{bestgreen}4.11e-03 ± 6.70e-04
& \cellcolor{bestgreen}1.92e-03 ± 2.80e-04
& \cellcolor{bestgreen}2.05e-03 ± 3.10e-04
& \cellcolor{bestgreen}1.61e-02 ± 4.20e-05
& \cellcolor{bestgreen}3.86e-03 ± 9.50e-04 \\

Improvement
& $\uparrow$ 11.58\%
& $\uparrow$ 89.04\%
& $\uparrow$ 81.18\%
& $\uparrow$ 48.88\%
& $\uparrow$ 79.06\%
& $\uparrow$ 36.30\% \\
\bottomrule
\end{tabular}
}
\end{table*}

\begin{figure*}[t]
\centering
\includegraphics[width=\textwidth]{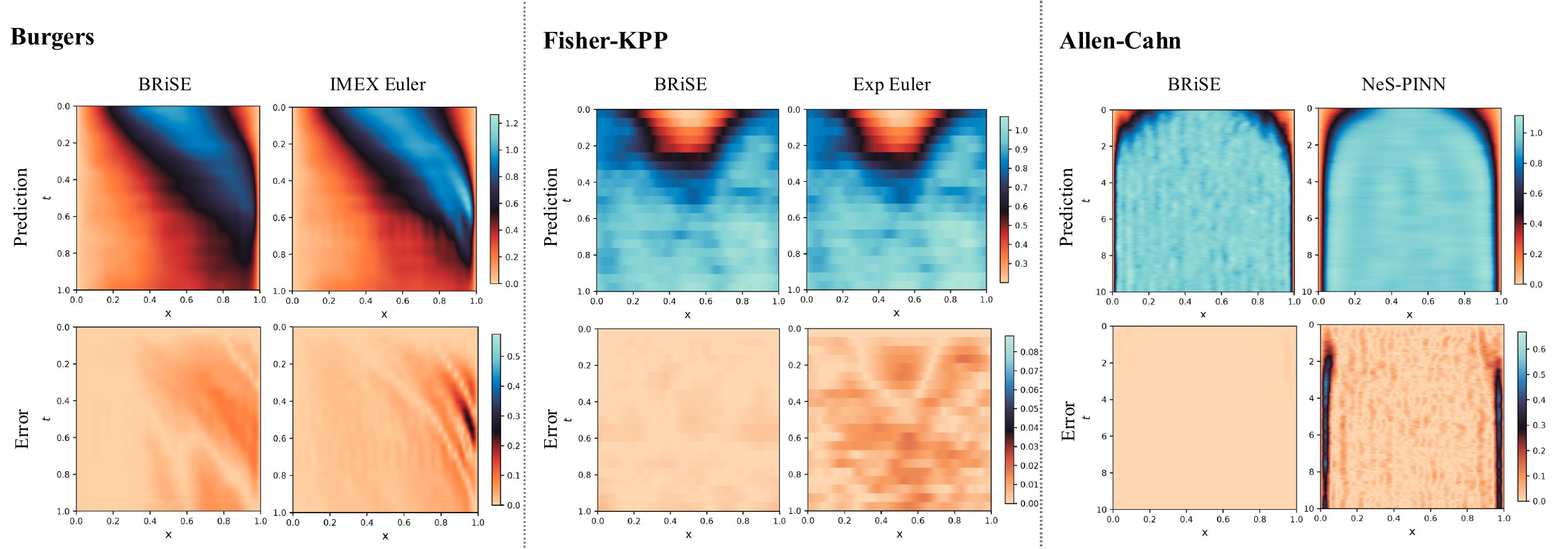}
\caption{Space-time predictions and error fields for the Burgers, Fisher--KPP, and Allen--Cahn benchmarks. For each benchmark, BRiSE is compared with the baseline indicated above the corresponding panels. The top row shows predicted solutions, and the bottom row shows their error fields.}
\label{fig:main_visualization}
\end{figure*}

\begin{figure}[htbp]
\centering
\includegraphics[width=\textwidth]{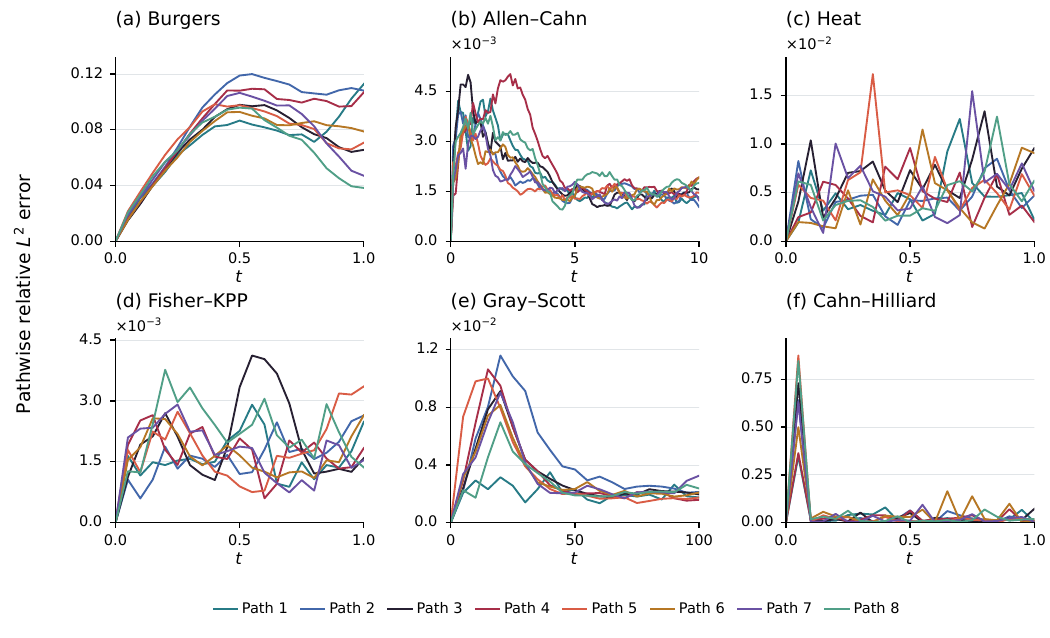}
\caption{Time evolution of the relative spatial $L^2$ error of BRiSE across six SPDE benchmarks. Each curve represents one of eight held-out Wiener realizations, with errors evaluated against the corresponding numerical reference solutions at each time point.}
\label{fig:pathwise_visualization}
\end{figure}

\paragraph{Main results.}
Table~\ref{tab:main_results} reports the space-time relative $L^2$ errors across six SPDE benchmarks. BRiSE achieves the lowest mean error on all six problems, consistently outperforming both physics-informed neural baselines and classical numerical solvers. Compared with the strongest baseline on each benchmark, BRiSE reduces the relative error by 89.04\% on Heat, 81.18\% on Fisher--KPP, and 79.06\% on Cahn--Hilliard. Substantial improvements are also observed on Allen--Cahn (48.88\%) and Gray--Scott (36.30\%), while BRiSE achieves a more modest gain of 11.58\% on Burgers. Among the baselines, Exponential Euler is generally the most competitive, whereas Euler--Maruyama fails to converge on four of the six benchmarks. These results demonstrate that BRiSE delivers consistently accurate pathwise approximations across diverse stochastic dynamics without requiring reference solution data during training. Besides, Figure~\ref{fig:main_visualization} visualizes the numerical solution between several numerical methods across different SPDE benchmarks. Moreover, Appendix~\ref{app:com_cost} discusses the computational cost of BRiSE is comparable to that of classical solvers, as it requires only a single network evaluation per coarse interval, while achieving superior accuracy. Figure~\ref{fig:pathwise_visualization} reports the temporal evolution of the relative spatial $L^2$ error across eight held-out Wiener realizations for each of the six benchmarks, complementing the aggregate space-time errors in Table~\ref{tab:main_results}.

\raggedbottom
\subsection{SPDE Benchmarks}\label{sec:SPDE}
We consider the following six SPDE benchmarks. Boundary conditions, spectral discretizations, and noise configurations are detailed in Appendix~\ref{app:reproducibility}.

\subsubsection{Example 1: Stochastic Burgers Equation.}
\begin{figure}[H]
\centering
\includegraphics[width=\textwidth]{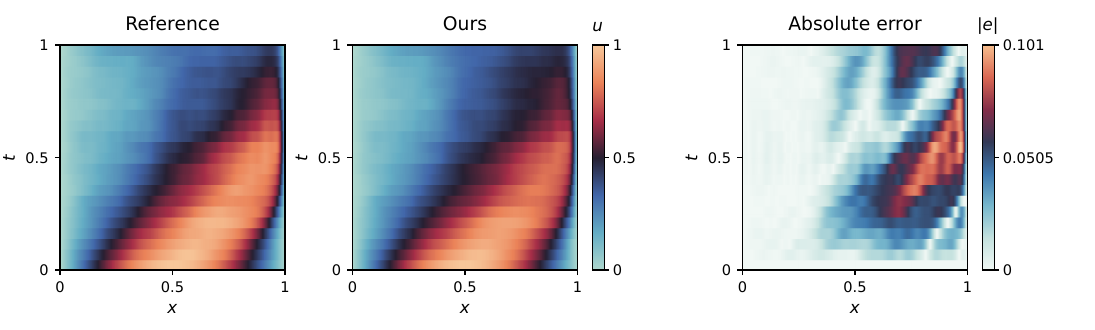}
\caption{Space-time visualization of the stochastic Burgers equation, comparing the numerical reference solution, BRiSE prediction, and absolute error for the same Wiener realization.}
\label{fig:burgers_visualization}
\end{figure}

We first consider the one-dimensional stochastic Burgers equation with nonlinear convection and additive noise:
\begin{equation}
\begin{aligned}
&\mathrm{d}u=\left(0.01u_{xx}-uu_x\right)\,\mathrm{d}t+0.1\,\mathrm{d}W_t^Q, \\
& u(0,x)=\sin(\pi x), \quad x\in(0,1),\quad t\in[0,1].
\end{aligned}
\end{equation}
The stochastic Burgers equation is a simplified model for nonlinear fluid transport subject to random external forcing. The nonlinear convection term describes wave steepening, while the diffusion term accounts for viscous dissipation. It provides a representative test of the interplay between nonlinear transport, diffusion, and stochastic perturbations. Figure~\ref{fig:burgers_visualization} compares the numerical reference solution and BRiSE prediction under the same Wiener realization, together with the pointwise absolute error.

\subsubsection{Example 2: Stochastic Allen--Cahn Equation}
\begin{figure}[H]
\centering
\includegraphics[width=\textwidth]{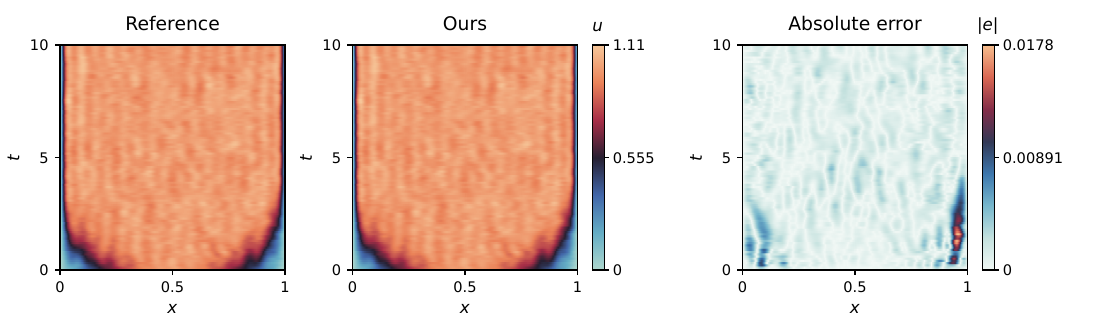}
\caption{Space-time visualization of the stochastic Allen--Cahn equation, comparing the numerical reference solution, BRiSE prediction, and absolute error for the same Wiener realization.}
\label{fig:ac_visualization}
\end{figure}

We next consider the stochastic Allen--Cahn equation, characterized by a cubic reaction term and weak diffusion:
\begin{equation}
\begin{aligned}
&\mathrm{d}u=\left(10^{-4}u_{xx}+u-u^3\right)\,\mathrm{d}t+0.01\,\mathrm{d}W_t^Q,\\
& u(0,x)=\sin(\pi x), \quad x\in(0,1),\quad t\in[0,10].
\end{aligned}
\end{equation}
The stochastic Allen--Cahn equation arises in the modeling of phase transitions and interface evolution in materials science. The cubic reaction term is associated with a double-well potential, while the diffusion term controls interfacial smoothing. Stochastic forcing represents random fluctuations that may influence phase evolution and interface dynamics. Figure~\ref{fig:ac_visualization} compares the numerical reference solution and BRiSE prediction under the same Wiener realization, together with the pointwise absolute error.

\subsubsection{Example 3: Stochastic Heat Equation}
\begin{figure}[H]
\centering
\includegraphics[width=\textwidth]{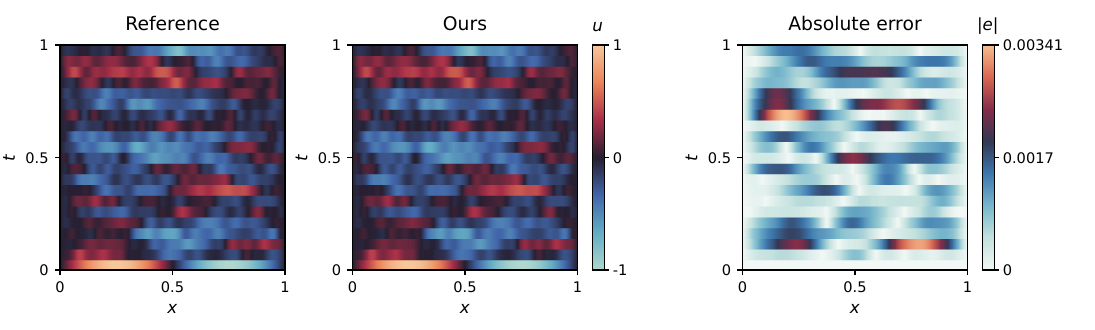}
\caption{Space-time visualization of the stochastic Heat equation, comparing the numerical reference solution, BRiSE prediction, and absolute error for the same Wiener realization.}
\label{fig:heat_visualization}
\end{figure}

We also consider a semilinear stochastic heat equation with nonlinear reaction and additive noise:
\begin{equation}
\begin{aligned}
&\mathrm{d}u=\left(u_{xx}-u-\sin u\right)\,\mathrm{d}t+\mathrm{d}W_t^Q,\\
& u(0,x)=\sin(2\pi x), \quad x\in(0,1),\quad t\in[0,1].
\end{aligned}
\end{equation}
The stochastic heat equation describes diffusion processes subject to random forcing. In this semilinear variant, the Laplacian governs spatial diffusion, while the nonlinear reaction term introduces state-dependent dissipation. The stochastic perturbation models random environmental fluctuations or external sources. Figure~\ref{fig:heat_visualization} compares the numerical reference solution and BRiSE prediction under the same Wiener realization, together with the pointwise absolute error.

\subsubsection{Example 4: Stochastic Fisher--KPP Equation}
\begin{figure}[htbp]
\centering
\includegraphics[width=\textwidth]{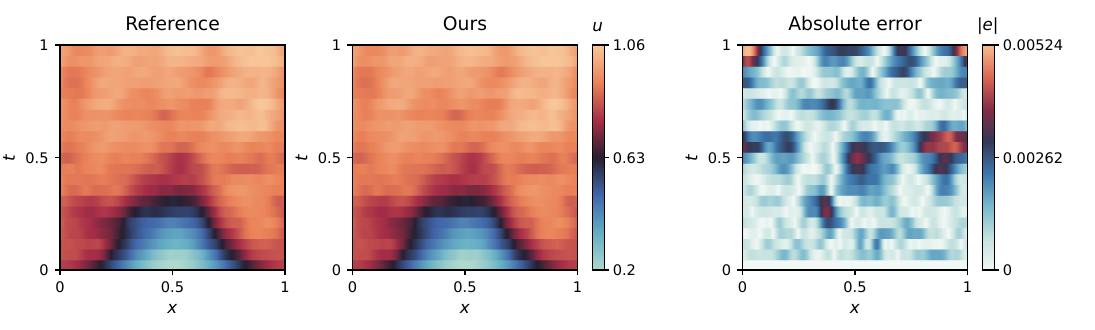}
\caption{Space-time visualization of the stochastic Fisher--KPP equation, comparing the numerical reference solution, BRiSE prediction, and absolute error for the same Wiener realization.}
\label{fig:fisher_kpp_visualization}
\end{figure}

To examine reaction--diffusion dynamics with nonlinear growth, we consider the stochastic Fisher--KPP equation:
\begin{equation}
\begin{aligned}
&\mathrm{d}u=\left(0.01u_{xx}+6u(1-u)\right)\,\mathrm{d}t+0.1\,\mathrm{d}W_t^Q,\\
&u(0,x)=0.5+0.3\cos(2\pi x), \quad x\in(0,1),\quad t\in[0,1].
\end{aligned}
\end{equation}
The stochastic Fisher--KPP equation is widely used to describe population growth and spatial spreading in random environments. The logistic reaction term models population growth with limited resources, while diffusion accounts for spatial dispersal. Random forcing captures environmental variability affecting population dynamics and propagating fronts. Figure~\ref{fig:fisher_kpp_visualization} compares the numerical reference solution and BRiSE prediction under the same Wiener realization, together with the pointwise absolute error.

\subsubsection{Example 5: Stochastic Gray--Scott Equation}
\begin{figure}[htbp]
\centering
\includegraphics[width=\textwidth]{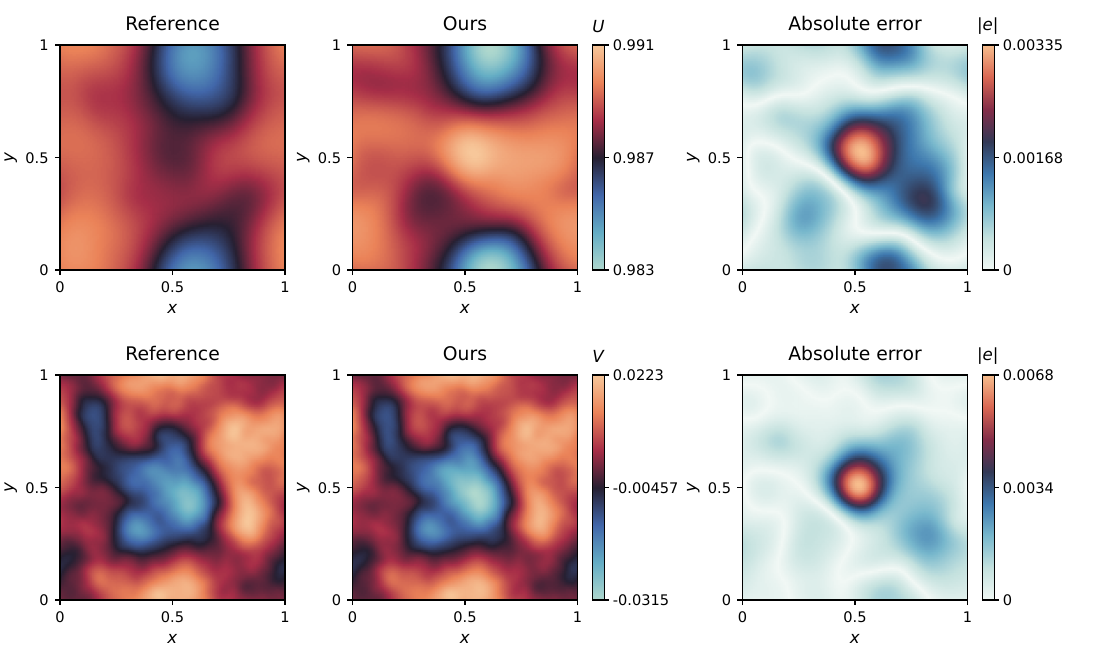}
\caption{Spatial snapshots of the 2D stochastic Gray--Scott system at $t=100$, comparing the numerical reference solution, BRiSE prediction, and pointwise absolute error for both $U$ (top row) and $V$ (bottom row) under the same Wiener realization.}
\label{fig:gray_scott_visualization}
\end{figure}

To investigate coupled nonlinear dynamics in two spatial dimensions, we consider the stochastic Gray--Scott system:
\begin{equation}
\begin{aligned}
&\mathrm{d}U=\left(0.0016\Delta U-UV^2+0.029(1-U)\right)\,\mathrm{d}t,\\
&\mathrm{d}V=\left(0.0008\Delta V+UV^2-0.086V\right)\,\mathrm{d}t +0.005\,\mathrm{d}W_t^{Q_V}. \\
&U(0,x,y)=1-0.5g(x,y), \quad V(0,x,y)=0.25g(x,y), \\
&(x,y)\in[0,1]^2,\quad t\in[0,100],
\end{aligned}
\end{equation}
where
\begin{equation}
g(x,y)=\exp\left( -\frac{d_{\mathbb T}(x,0.5)^2+d_{\mathbb T}(y,0.5)^2}{2(0.1)^2} \right), \quad d_{\mathbb T}(a,b)=\min\{|a-b|,\,1-|a-b|\}.
\end{equation}
Here $d_{\mathbb T}(a,b)$ is the shortest distance between $a$ and $b$ on the unit interval with its endpoints identified, so $g$ uses the periodic distance in each spatial coordinate.
The stochastic Gray--Scott system models autocatalytic reactions on the product space $\mathcal H=L^2(D;\mathbb R^2)$, with state $u=(U,V)$. Its linear part is $A=\operatorname{diag}(0.0016\Delta,0.0008\Delta)$, and its displayed noise can be written with $\sigma=0.005$ and covariance $\operatorname{diag}(0,Q_V)$ on the product space, where $Q_V$ is the covariance in the $V$ equation. Thus the absence of direct noise in $U$ is represented by a zero covariance block. The common-eigenbasis condition in Section~\ref{sec:preliminaries} still requires compatible Laplacian and covariance eigenfunctions. The system describes chemical reactions coupled with spatial diffusion. The variables $U$ and $V$ represent the concentrations of two interacting chemical species, whose nonlinear reaction kinetics can generate complex spatial patterns. The stochastic forcing accounts for random fluctuations in the reaction--diffusion dynamics. Figure~\ref{fig:gray_scott_visualization} compares the numerical reference solution and BRiSE prediction under the same Wiener realization, together with the pointwise absolute error.

\subsubsection{Example 6: Stochastic Cahn--Hilliard Equation}
\begin{figure}[htbp]
\centering
\includegraphics[width=\textwidth]{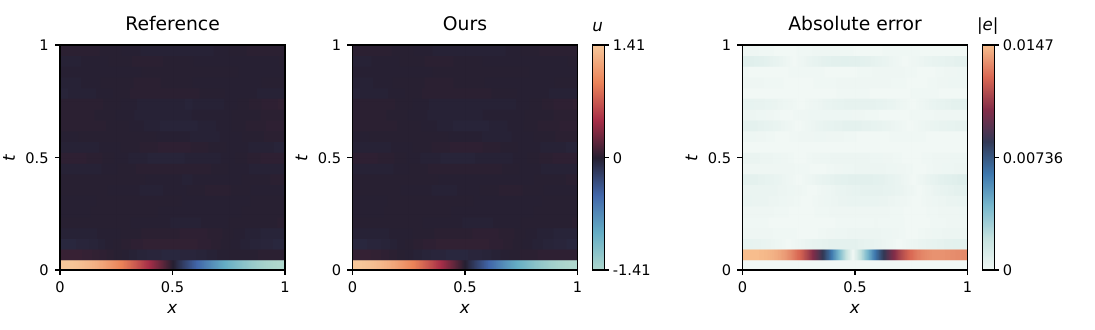}
\caption{Space-time visualization of the stochastic Cahn--Hilliard equation, comparing the numerical reference solution, BRiSE prediction, and absolute error for the same Wiener realization.}
\label{fig:cahn_hilliard_visualization}
\end{figure}

Finally, we consider the stochastic Cahn--Hilliard equation, which involves a fourth-order spatial differential operator:
\begin{equation}
\begin{aligned}
&\mathrm{d}u=\left(-u_{xxxx}+\partial_{xx}(u^3-u)\right)\,\mathrm{d}t+\mathrm{d}W_t^Q,\\
&u(0,x)=\sqrt{2}\cos(\pi x), \quad x\in(0,1),\quad t\in[0,1].
\end{aligned}
\end{equation}
The stochastic Cahn--Hilliard equation arises in the description of phase separation and coarsening phenomena in binary mixtures. The fourth-order diffusion operator and nonlinear chemical potential govern the evolution of interfaces between different phases. Stochastic forcing represents random fluctuations affecting the phase-separation process.
Figure~\ref{fig:cahn_hilliard_visualization} compares the numerical reference solution and BRiSE prediction under the same Wiener realization, together with the pointwise absolute error.

\section{Additional Experiments}
\label{sec:ab}
In this section, we conduct hyperparameter analysis and ablation studies on the stochastic Fisher-KPP equation to investigate the sensitivity of BRiSE to its model configurations and assess the contributions of its key components. Specifically, we examine the effects of temporal resolution and network width, followed by ablations on correction order, OU bridge integral features, and stochastic convolution factorization. All variants are evaluated on the same eight held-out Wiener paths using a common space-time evaluation grid. We further assess robustness of our method to noise intensity and spatial regularity on the stochastic Burgers equation.

\subsection{Hyperparameter Analysis}

\subsubsection{Time interval Sensitivity}
We investigate the effect of temporal resolution by varying the number of time intervals $N\in\{10,20,40,80\}$. Figure~\ref{fig:ablation_visualization} (a) shows that as $N$ increases, the space-time relative $L^2$ error consistently decreases from $4.28\times10^{-3}$ to $4.94\times10^{-4}$, corresponding to an 88.5\% reduction. These results demonstrate that finer temporal discretization substantially improves the accuracy of BRiSE in approximating stochastic dynamics.

\subsubsection{Network Capacity Sensitivity}
We investigate the effect of network capacity by varying the hidden width $d\in\{64,128,256,512\}$. Figure~\ref{fig:ablation_visualization} (b) shows that the space-time relative $L^2$ errors remain within a narrow range of $1.85\times10^{-3}$ to $2.09\times10^{-3}$, with no consistent improvement as the network width increases. Notably, increasing $d$ from 128 to 512 yields only a marginal reduction in error, from $1.86\times10^{-3}$ to $1.85\times10^{-3}$. These results suggest that BRiSE achieves competitive accuracy with relatively compact networks and exhibits limited sensitivity to network width under the considered settings.

\subsection{Ablation Study}
\begin{figure*}[t]
\centering
\includegraphics[width=\textwidth]{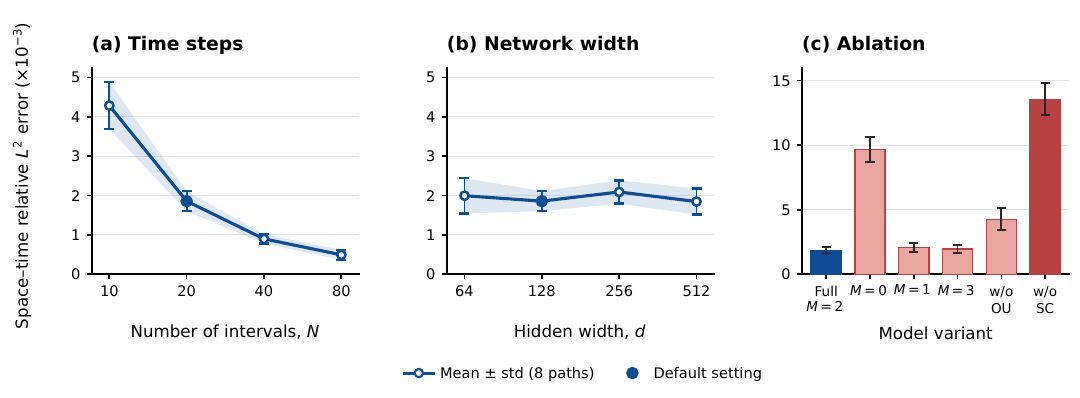}
\caption{Hyperparameter sensitivity and ablation studies on the stochastic Fisher-KPP equation. (a) Effect of the number of time intervals $N$. (b) Effect of network width $d$. (c) Ablation of correction order $M$, OU bridge integral features (OU), and stochastic convolution factorization (SC). All results are reported as mean ± std of the space-time relative $L^2$ error over eight held-out Wiener paths. Blue markers indicate the default configuration ($N=20$, $d=128$, $M=2$).}
\label{fig:ablation_visualization}
\end{figure*}

We conduct ablation studies on the stochastic Fisher-KPP equation to evaluate the contributions of the key components of BRiSE, and the results are shown in Figure~\ref{fig:ablation_visualization} (c). With $M=0$, the decoder retains the semigroup propagation and the frozen drift term $D_sG_n$, and the reconstruction adds the stochastic convolution. This increases the space-time relative $L^2$ error from $1.86\times10^{-3}$ to $9.65\times10^{-3}$. Increasing the correction order to $M=2$ yields the lowest error, while $M=3$ provides no further improvement. Removing the integrated innovation feature $\widehat J_n$ while retaining the endpoint innovation $\eta_n$ increases the error to $4.25\times10^{-3}$, highlighting the importance of capturing intra-interval stochastic information beyond endpoint innovations. Removing the explicit decomposition of the solution into the remainder and stochastic convolution leads to a substantial error increase to $1.35\times10^{-2}$, approximately 7.26 times that of the full model. These results demonstrate the effectiveness of the proposed stochastic feature representation and factorized solution structure.

\subsection{Robustness to Noise Intensity and Spatial Regularity}\label{sec:noise_robustness}

To assess the robustness of BRiSE under varying stochastic forcing conditions, we investigate its sensitivity to noise intensity and spatial regularity using the stochastic Burgers equation as below,
\begin{equation}
\begin{aligned}
&\mathrm{d}u =\left(0.01\,u_{xx}-uu_x\right)\,\mathrm{d}t +\sigma\,\mathrm{d}W_t^{\widetilde Q_{\alpha_Q},32},\\
&u(0,x)=\sin(\pi x),\quad x\in[0,1],
\end{aligned}
\end{equation}
with homogeneous Dirichlet boundary conditions and $T=1$. The truncated $Q$-Wiener process is defsined as
\begin{equation}
W_t^{\widetilde Q_{\alpha_Q},32}(x)= \sum_{j=1}^{32} \sqrt{\widetilde q_j^{(\alpha_Q)}}\, \beta_j(t)e_j(x),
\end{equation}
where $e_j(x)=\sqrt{2}\sin(j\pi x)$ and the covariance eigenvalues are normalized as
\begin{equation}
\widetilde q_j^{(\alpha_Q)}
  = \frac{S(1)}{S(\alpha_Q)}\,(1+j^2)^{-\alpha_Q},
\quad
S(\alpha) = \sum_{k=1}^{32}(1+k^2)^{-\alpha}.
\end{equation}
Here, $\sigma$ controls the noise intensity, while $\alpha_Q$ governs its spatial regularity through the decay of the covariance eigenvalues, with larger $\alpha_Q$ corresponding to smoother spatial noise. The covariance $\widetilde Q_{\alpha_Q}$ has eigenvalues $\widetilde q_j^{(\alpha_Q)}$ on the 32 retained sine modes and zero on all other modes. The normalization ensures that $\operatorname{Tr}(\widetilde Q_{\alpha_Q})$ remains constant across different values of $\alpha_Q$, allowing us to isolate the effect of spatial regularity from changes in total noise variance. We investigate these two factors separately by varying $\sigma$ at fixed $\alpha_Q=1$ and varying $\alpha_Q$ at fixed $\sigma=0.1$.

\paragraph{Experimental setup.}
We conduct two controlled experiments: (i) varying the noise intensity $\sigma\in\{0.05,0.1,0.2,0.3\}$ with fixed $\alpha_Q=1$, and (ii) varying the spatial regularity $\alpha_Q\in\{0,1,2\}$ with fixed $\sigma=0.1$. All configurations employ 32 sine modes, 256 spatial quadrature points, and $N=20$ time intervals. BRiSE uses two hidden layers of width 128 with tanh activations and correction order $M=2$. We compare BRiSE against 20-step Exponential Euler on eight held-out Wiener paths. Results are reported as the mean and sample standard deviation of the space--time relative $L^2$ error.

\begin{table}[htbp]
\centering
\caption{Robustness to noise intensity and spatial regularity on the stochastic Burgers equation.
Results are reported as mean ± std of the space--time relative $L^2$ error over eight held-out Wiener paths; lower is better. Better results are highlighted in \colorbox{bestgreen}{green}, while failures (errors exceeding 1) are highlighted in \colorbox{failred}{red}.}
\label{tab:noise_robustness}
\begin{tabular}{lcc}
\toprule
Noise setting & BRiSE & Exp Euler \\
\midrule
\multicolumn{3}{l}{\textit{Noise intensity ($\alpha_Q=1$)}} \\
$\sigma=0.05$
& \cellcolor{bestgreen}{3.94e-02 ± 2.37e-03}
& 1.26e-01 ± 2.42e-02 \\
$\sigma=0.10$
& \cellcolor{bestgreen}{7.54e-02 ± 8.50e-03}
& 2.20e-01 ± 1.42e-01 \\
$\sigma=0.20$
& \cellcolor{failred}$>1$
& \cellcolor{failred}$>1$ \\
$\sigma=0.30$
& \cellcolor{failred}$>1$
& \cellcolor{failred}$>1$ \\
\midrule
\multicolumn{3}{l}{\textit{Noise regularity ($\sigma=0.10$)}} \\
$\alpha_Q=0$
& \cellcolor{bestgreen}{7.17e-02 ± 3.61e-03}
& \cellcolor{failred}$>1$ \\
$\alpha_Q=1$
& \cellcolor{bestgreen}{7.54e-02 ± 8.50e-03}
& 2.20e-01 ± 1.42e-01 \\
$\alpha_Q=2$
& \cellcolor{bestgreen}{4.65e-02 ± 5.42e-03}
& 9.13e-02 ± 1.97e-02 \\
\bottomrule
\end{tabular}
\end{table}

\subsubsection{Noise Intensity}
We first investigate the effect of noise intensity by varying $\sigma\in\{0.05,0.1,0.2,0.3\}$ while fixing $\alpha_Q=1$. Table~\ref{tab:noise_robustness} shows that as $\sigma$ increases from $0.05$ to $0.1$, the space--time relative $L^2$ error of BRiSE increases from $3.94\times10^{-2}$ to $7.54\times10^{-2}$, reflecting the growing difficulty of resolving stronger stochastic fluctuations and their interactions with nonlinear dynamics. Nevertheless, BRiSE reduces the error by approximately 69\% and 66\% compared with Exponential Euler at these two noise levels, respectively. Under stronger stochastic forcing ($\sigma\geq0.2$), BRiSE fails to complete training, while Exponential Euler exhibits errors exceeding one. These results demonstrate the accuracy advantage of BRiSE under moderate noise intensities, while revealing the limitations of its training stability as stochastic forcing becomes stronger.

\subsubsection{Noise Regularity}
We further investigate the effect of spatial noise regularity by varying $\alpha_Q\in\{0,1,2\}$ while keeping $\sigma=0.1$ and the covariance trace fixed. Unlike noise intensity, spatial regularity controls how stochastic energy is distributed across spatial frequencies. Smaller values of $\alpha_Q$ lead to slower spectral decay and stronger high-frequency fluctuations, making the stochastic forcing spatially rougher. As shown in Table~\ref{tab:noise_robustness}, at $\alpha_Q=0$, Exponential Euler produces errors exceeding one, whereas BRiSE maintains a space--time relative $L^2$ error of $7.17\times10^{-2}$. As the noise becomes smoother at $\alpha_Q=2$, BRiSE achieves a lower error of $4.65\times10^{-2}$, compared with $9.13\times10^{-2}$ for Exponential Euler. Notably, BRiSE exhibits comparable accuracy at $\alpha_Q=0$ and $\alpha_Q=1$, despite their substantially different spectral distributions. These findings suggest that BRiSE can accommodate changes in spatial noise regularity while maintaining accurate pathwise approximations, even under rougher stochastic forcing.

\section{Discussions and Conclusions}
\label{sec:conclusion}
\paragraph{What does BRiSE achieve?}
In this work, we introduce BRiSE, a structure-aware neural time integrator for semilinear SPDEs with additive noise. Across six SPDE benchmarks, BRiSE consistently achieved the lowest mean space-time relative $L^2$ errors, reducing errors by 11.6\%--89.0\% compared with the strongest competing methods. Notably, BRiSE requires no reference solution trajectories for training and generalizes to unseen stochastic realizations without path-specific optimization. With only one neural network evaluation per time interval, BRiSE provides an accurate and reusable framework for pathwise stochastic simulation.

\paragraph{Why does BRiSE work?}
The success of BRiSE suggests a broader insight: accurate stochastic simulation may benefit more from learning the right dynamics than from learning the entire solution. By explicitly separating the stochastic convolution, BRiSE shifts the learning burden toward nonlinear corrections that depend on individual noise realizations. Yet this separation alone is insufficient, as the nonlinear evolution remains sensitive to stochastic fluctuations within each time interval. OU bridge-informed features address this challenge by preserving information beyond endpoint innovations, while the semigroup-factorized decoder translates the learned corrections into continuous-time trajectories. Our ablation studies reinforce this perspective, showing substantial accuracy degradation when either stochastic convolution factorization or bridge information is removed. These findings highlight the importance of combining analytical structure with selective neural approximation, rather than treating stochastic dynamics as an entirely black-box learning problem.

\paragraph{Where can BRiSE go next?}
Looking ahead, the principles underlying BRiSE offer opportunities for extending structure-aware neural solvers to broader classes of stochastic systems. The current framework focuses on semilinear SPDEs with additive noise and tractable spectral representations, where stochastic convolutions can be handled explicitly. Extending this approach to multiplicative noise, more general differential operators, and complex spatial geometries presents promising directions for future research. More broadly, BRiSE illustrates a perspective for scientific machine learning in which known mathematical structures guide neural approximation, shifting the emphasis from learning entire stochastic solutions to learning the dynamics that genuinely require approximation.

\section*{Data Availability}
The datasets generated during and/or analysed during the current study are available from the corresponding author on request.

\section*{Conflict of interest}
The authors declared that they have no conflict of interest.

\section*{Acknowledgments}
This work was supported by the National Natural Science Foundation of China (12671434, 12571391, 12271366).

\bibliography{references}
\bibliographystyle{plainnat}

\newpage
\appendix
\onecolumn
\flushbottom
\addtocontents{toc}{\protect\setcounter{tocdepth}{2}}

{\centering\Large\bfseries
Learn the Dynamics, Not the Noise:\\ A Semigroup-Factorized Solver for Stochastic PDEs -- Appendix\par}
\label{sec:appendix}

\noindent\rule{\textwidth}{0.6pt}
\begingroup
\renewcommand{\contentsname}{{\Large Contents}}
\let\appendixtocaddvspace\addvspace
\renewcommand{\addvspace}[1]{\appendixtocaddvspace{0.65em}}
\tableofcontents
\endgroup
\noindent\rule{\textwidth}{0.6pt}

\section{Additional Mathematical Details}
\label{app:math}
In this section, we provide supplementary mathematical details of BRiSE, including the construction of OU path features, the analytical derivation of semigroup-integrated temporal basis functions, the structural consistency of the neural decoder, and the numerical evaluation of pathwise mild-equation residuals.

\subsection{OU Path Construction} \label{app:OU Path Construction}
We start from the modal OU dynamics in Equation~\eqref{eq:modal_ou},
\begin{equation}\label{eq:app_modal_ou}
\mathrm{d}z_k(t)=-\lambda_k z_k(t)\,\mathrm{d}t +\sigma\sqrt{q_k}\,d\beta_k(t),\quad z_k(0)=0,
\end{equation}
where $\{\beta_k\}_{k=1}^{K}$ are independent standard Brownian motions. Applying the integrating factor $e^{\lambda_k t}$ to Equation~\eqref{eq:app_modal_ou} over $[t,t+h]$ gives the exact transition
\begin{equation}
z_k(t+h)=e^{-\lambda_k h}z_k(t)+\sigma\sqrt{q_k}\,b_k(h)\xi_k,
\end{equation}
where $\xi_k\sim\mathcal N(0,1)$. By the It\^o isometry, $b_k(h)^2=\int_0^h e^{-2\lambda_k r}\,\mathrm{d}r$, so
\begin{equation}
b_k(h)= \begin{cases}
\displaystyle \sqrt{\frac{1-e^{-2\lambda_k h}}{2\lambda_k}}, & \lambda_k>0,\\
\sqrt{h}, & \lambda_k=0.
\end{cases}
\end{equation}
The Gaussian innovations are sampled independently across spectral modes and disjoint time increments. This transition enables exact sampling of the truncated OU process at discrete time points.

For each coarse interval $[t_n,t_n+H]$, subtracting the propagated initial state as in Equation~\eqref{eq:innovation} leaves the stochastic integral
\begin{equation}
Y_{n,k}(s)= \sigma\sqrt{q_k} \int_{t_n}^{t_n+s} e^{-\lambda_k(t_n+s-r)}\,d\beta_k(r).
\end{equation}
To construct the integrated innovation feature, we sample the OU trajectory on a subgrid $0=\zeta_0<\cdots<\zeta_L=H$ and apply the trapezoidal rule to the integral defining $J_{n,k}$ in Equation~\eqref{eq:ou_features}, giving
\begin{equation}
\widehat J_{n,k} = \sum_{\ell=0}^{L-1} \frac{\zeta_{\ell+1}-\zeta_\ell}{2} \left( Y_{n,k}(\zeta_\ell)+Y_{n,k}(\zeta_{\ell+1}) \right).
\end{equation}
While the OU transitions are exact at the sampled time points, the integrated innovation remains a numerical approximation. Importantly, the path features and the mild-equation residual are evaluated using the same OU realization, ensuring consistent stochastic forcing throughout training.

\subsection{Semigroup-Integrated Temporal Basis.} \label{app:temporal_basis}
Substituting $p_j(r/H)=(r/H)^j$ into Equation~\eqref{eq:basis_integrals} and integrating by parts for $j=1,2$ gives, for $\lambda_k>0$,
\begin{equation}
\begin{aligned}
& \Psi_{1,k}(s;H) = \frac{\lambda_k s-1+e^{-\lambda_k s}} {H\lambda_k^2}, \\
& \Psi_{2,k}(s;H)=\frac{\lambda_k^2s^2-2\lambda_ks+ 2(1-e^{-\lambda_ks})}{H^2\lambda_k^3}.
\end{aligned}
\end{equation}
When $\lambda_k=0$, the kernel is one, so integrating the two monomials directly gives
\begin{equation}\label{eq:zero_mode_basis}
\Psi_{1,k}(s;H)=\frac{s^2}{2H}, \quad \Psi_{2,k}(s;H)=\frac{s^3}{3H^2}.
\end{equation}
To improve numerical stability for small $\lambda_k s$, the basis functions can equivalently be expressed using the exponential integrator $\varphi$-functions
\begin{equation}
\Psi_{j,k}(s;H) = \frac{j!\,s^{j+1}}{H^j} \varphi_{j+1}(-\lambda_ks),
\end{equation}
where
\begin{equation}
\varphi_m(x)=\sum_{\ell=0}^{\infty}\frac{x^\ell}{(\ell+m)!}.
\end{equation}
These analytical basis functions depend only on the linear operator, the interval length, and the evaluation time, allowing precomputation on fixed temporal grids.

\paragraph{Structural consistency.}
Evaluating the decoder in Equation~\eqref{eq:brise_decoder} and its derivative at $s=0$ gives two consistency properties. Since $E_0=I$, $D_0=0$, and $\Psi_j(0;H)=0$, the reconstructed trajectory satisfies the initial-state consistency condition $\widehat a_{\theta,n}(0)=x_n$. Moreover, the temporal basis functions satisfy $\partial_s\Psi_j(0;H)=0$ for $j\geq1$, ensuring initial-drift consistency:
\begin{equation}
\left.\frac{\mathrm{d}}{\mathrm{d}s}\widehat a_{\theta,n}(s)\right|_{s=0}=-\Lambda x_n+G_n.
\end{equation}
These properties hold independently of the learned coefficients, ensuring that neural corrections preserve the prescribed initial state and the known initial drift.

\subsection{Numerical Residual Evaluation and Rollout Training} \label{app:numerical_residual_evaluation}
The residual uses exponential product integration with a piecewise-linear interpolant of the nonlinear forcing. For a collocation time $s_m$, let $0=r_{m,0}<\cdots<r_{m,Q_m}=s_m$ be the partition used for the mild integral, and set $g_{\theta,n}(r)=G(\widehat a_{\theta,n}(r),z(t_n+r))$. Let $\ell_{m,i}$ be the piecewise-linear basis function that is one at $r_{m,i}$ and zero at the other nodes. Substituting $\sum_i \ell_{m,i}(r)g_{\theta,n}(r_{m,i})$ for the forcing in Equation~\eqref{eq:brise_residual} and integrating each basis function against the semigroup gives
\begin{equation}\label{eq:app_residual_quadrature}
\begin{aligned}
\widehat{\mathcal R}_{\theta,n}(s_m)
&=\widehat a_{\theta,n}(s_m)-E_{s_m}x_n
   -\sum_{i=0}^{Q_m} W_{m,i}\,g_{\theta,n}(r_{m,i}),\\
W_{m,i}&=\int_0^{s_m} E_{s_m-r}\ell_{m,i}(r)\,\mathrm{d}r .
\end{aligned}
\end{equation}
The diagonal matrices $W_{m,i}$ account for the exponential kernel in each mode, including the zero mode if present. All collocation and quadrature states are reconstructed from the same interval-wise coefficient prediction. The OU values are taken from the realization used to construct the path features. With approximate path information, $g_{\theta,n}$ uses the same $\widetilde z$ as the decoder and its inputs.

Training initializes the projected remainder state from the prescribed initial condition and recursively propagates the predicted interval endpoints. Gradients are backpropagated through the unrolled trajectory and nonlinear evaluations while holding each sampled noise realization fixed. Fresh noise paths are generated across optimization iterations, and model selection is based on held-out full-horizon physics residuals rather than reference solution errors.

\section{Experimental and Reproducibility Details}
\label{app:reproducibility}
In this section, we provide additional implementation details of BRiSE, including the spectral discretization, neural network architecture, and training configurations used across the six SPDE benchmarks to facilitate reproducibility.

\subsection{Spectral Discretization}
BRiSE employs an orthonormal spectral Galerkin discretization, where the solution states are represented by real-valued spectral coefficients. The Galerkin linear operator in coefficient coordinates is $A_K=I_K^{-1}P_K A I_K=-\Lambda$, allowing its semigroup $E_s=e^{-\Lambda s}$ to be evaluated analytically. Nonlinear terms are computed by reconstructing the physical fields on spatial quadrature grids and projecting the resulting nonlinearities onto the retained spectral modes. The spectral bases are chosen according to the boundary conditions of each SPDE, including Fourier sine, cosine, and periodic Fourier bases. For two-dimensional Gray-Scott, complex Fourier transforms are used internally, while the neural network operates on real-valued spectral coordinates. The spatial discretization and truncated noise configurations are summarized in Table~\ref{tab:spectral_config}. Then we will introduce different spectral basis we used in the experiments.

\paragraph{Fourier sine basis.}
For Burgers, Allen--Cahn, and Heat, we employ the orthonormal Fourier sine basis on $D=(0,1)$,
\begin{equation}
\phi_k(x)=\sqrt{2}\sin(k\pi x), \quad k=1,\ldots,32.
\end{equation}
This basis satisfies the homogeneous Dirichlet boundary conditions. For the diffusion operator $A=\nu\partial_{xx}$, the corresponding eigenvalues are $\lambda_k=\nu(k\pi)^2$, enabling diagonal semigroup propagation. The nonlinear terms are evaluated on spatial quadrature grids and projected onto the 32 retained modes.

\paragraph{Fourier cosine basis.}
For Cahn--Hilliard, we employ the orthonormal Fourier cosine basis on $D=(0,1)$,
\begin{equation}
\phi_k(x)=\sqrt{2}\cos(k\pi x), \quad k=1,\ldots,32.
\end{equation}
This basis satisfies the homogeneous boundary conditions $u_x=u_{xxx}=0$. The constant mode is excluded, restricting the spectral representation to zero-mean states. The fourth-order linear operator is diagonal in this basis, with eigenvalues $\lambda_k=(k\pi)^4$. Nonlinear terms are evaluated using midpoint quadrature and projected onto the retained cosine modes.

\paragraph{One-dimensional periodic Fourier basis.}
For Fisher--KPP, we employ the one-dimensional periodic Fourier basis on $D=(0,1)$,
\begin{equation}
\phi_0(x)=1,\quad \phi_k^c(x)=\sqrt{2}\cos(2\pi kx),\quad \phi_k^s(x)=\sqrt{2}\sin(2\pi kx),
\end{equation}
with $k=1,\ldots,16$. The resulting representation contains one constant mode and 16 sine-cosine pairs, yielding 33 real spectral coefficients. The constant mode is retained with zero linear decay rate, for which the semigroup-integrated basis functions are evaluated using Equation~\eqref{eq:zero_mode_basis}.

\paragraph{Two-dimensional periodic Fourier basis.}
For the Gray--Scott system, we adopt a two-dimensional periodic Fourier representation on $D=(0,1)^2$, retaining wavevectors $\boldsymbol{k}=(k_x,k_y)$ satisfying $$ |k_x|\leq15,\quad |k_y|\leq15. $$ The Fourier coefficients satisfy conjugate symmetry for real-valued fields and are packed into an equivalent real-valued representation. Each species contains $31^2=961$ real spectral degrees of freedom, resulting in a combined state dimension of 1,922. Spatial reconstruction and nonlinear projection are performed using fast Fourier transforms (FFTs) on a $64\times64$ grid.

\begin{table}[htbp]
\centering
\caption{Spectral discretization and noise configurations
of BRiSE across six SPDE benchmarks. Here, $K$ denotes the number of real-valued spectral coefficients, $q_k$ denotes the noise covariance eigenvalues, and $\sigma$ is the noise intensity.}
\label{tab:spectral_config}
\small \setlength{\tabcolsep}{5pt}
\begin{adjustbox}{max width=\textwidth}
\begin{tabular}{llccl}
\toprule
\textbf{SPDE}
& \textbf{Spectral basis}
& $\boldsymbol{K}$
& \textbf{Spatial grid}
& \textbf{Noise spectrum} \\
\midrule
Burgers
& Fourier sine
& 32
& 256
& $q_k=(1+k^2)^{-1}$\\

Allen--Cahn
& Fourier sine
& 32
& 128
& $q_k=1$\\

Heat
& Fourier sine
& 32
& 256
& $q_k=1$\\

Fisher--KPP
& Periodic Fourier
& 33
& 256
& $q_0=1,\quad q_k=(1+k^2)^{-1}$\\

Gray--Scott
& Periodic Fourier
& 1922
& $64\times64$
& $q_{\boldsymbol{k}}^{U}=0,\quad
q_{\boldsymbol{k}}^{V}=(1+|\boldsymbol{k}|^2)^{-2}$\\

Cahn--Hilliard
& Fourier cosine
& 32
& 256
& $q_1=0,\quad q_k=[k(\log k)^2]^{-1},\ k\geq2$\\
\bottomrule
\end{tabular}
\end{adjustbox}
\end{table}

\subsection{Network Architecture}
BRiSE employs a global feedforward multilayer perceptron (MLP) with Tanh activations to predict the temporal correction coefficients in the spectral domain. Rather than processing each spectral mode independently, the network jointly takes all retained spectral coefficients and stochastic path features as input. Specifically, for $K$ spectral modes, the input consists of the current remainder state $x_n$, stochastic convolution state $z_n$, OU endpoint innovation $\eta_n$, numerically integrated OU innovation $\widehat{J}_n$, and interval length $H$, resulting in an input dimension of $4K+1$. The network outputs $2K$ coefficients corresponding to the linear and quadratic temporal correction channels.

All hidden layers have width 128, with two hidden layers for Burgers, Allen--Cahn, Heat, and Fisher--KPP, and five for Gray--Scott and Cahn--Hilliard. The output layer uses a linear activation, with its weights and biases initialized to zero, while the remaining layers retain the default PyTorch initialization. This initialization yields zero temporal corrections before training, while the analytical semigroup terms remain active. The detailed network configurations are summarized in Table~\ref{tab:network_config}.

\begin{table}[htbp]
\centering
\caption{Network architecture configurations of BRiSE
across six SPDE benchmarks. }
\label{tab:network_config}
\small \setlength{\tabcolsep}{4pt}
\begin{tabular}{lcccc}
\toprule
\textbf{SPDE}
& \textbf{Input}
& \textbf{Depth}
& \textbf{Width}
& \textbf{Params} \\
\midrule
Burgers          & 129  & 2 & 128 & 41,408 \\
Allen--Cahn    & 129  & 2 & 128 & 41,408 \\
Heat   & 129  & 2 & 128 & 41,408 \\
Fisher--KPP      & 133  & 2 & 128 & 42,178 \\
Gray--Scott   & 7689 & 5 & 128 & 1,546,244 \\
Cahn--Hilliard   & 129  & 5 & 128 & 90,944 \\
\bottomrule
\end{tabular}
\end{table}

\subsection{Training Configuration}
We train a separate BRiSE model for each SPDE benchmark using the Adam optimizer \cite{kingma2014adam} and the pathwise mild-equation residual, without requiring reference solution trajectories. At each optimization iteration, fresh stochastic paths are sampled to construct the OU features and evaluate the corresponding physics residuals. Each coarse time interval is divided into 32 subintervals for stochastic path construction, with eight nonzero temporal collocation points used for residual evaluation. The sampled loss in Equation~\eqref{eq:mild_loss} uses the quadrature residual in Equation~\eqref{eq:app_residual_quadrature}. Training employs fully differentiable rollouts, where predicted endpoint states are propagated across consecutive intervals without detaching gradients.

For model selection, we evaluate the full-horizon physics residual on an independent set of validation paths every 50 optimization updates. Validation uses a refined temporal grid with 64 subdivisions and 64 nonzero collocation points per interval. The checkpoint achieving the lowest validation residual is selected, without accessing reference solutions. Most benchmarks are trained for 1,500 updates using a learning rate of $10^{-3}$ and a batch size of 16, while Gray--Scott employs a learning rate of $10^{-5}$ and a batch size of 4. For Burgers, we adopt a rollout curriculum that covers the first 10 intervals during the initial 200 updates before extending to the full 20-interval horizon; all other benchmarks use full-horizon training throughout. The problem-specific training configurations are summarized in Table~\ref{tab:training_config}.

\begin{table}[htbp]
\centering
\caption{Training configurations of BRiSE across six SPDE benchmarks.
$H$ denotes the coarse time step and $N$ the number of intervals.}
\label{tab:training_config}
\small \setlength{\tabcolsep}{8pt}
\begin{tabular}{lccccc}
\toprule
\textbf{SPDE}
& $\boldsymbol{H}$
& $\boldsymbol{N}$
& \textbf{Learning Rate}
& \textbf{Batch Size}
& \textbf{Iters} \\
\midrule
Burgers               & 0.05 & 20  & $10^{-3}$ & 16 & 1,500 \\
Allen--Cahn  & 0.1  & 100 & $10^{-3}$ & 16 & 1,500 \\
Heat         & 0.05 & 20  & $10^{-3}$ & 16 & 1,500 \\
Fisher--KPP            & 0.05 & 20  & $10^{-3}$ & 16 & 1,500 \\
Gray--Scott         & 5.0  & 20  & $10^{-5}$ & 4  & 1,500 \\
Cahn--Hilliard         & 0.05 & 20  & $10^{-3}$ & 16 & 1,500 \\
\bottomrule
\end{tabular}
\end{table}

\subsection{Data Generator}
To construct numerical reference solutions, we first project each SPDE onto a finite-dimensional spectral Galerkin space, yielding a system of stochastic evolution equations for the retained spectral coefficients. The driving Wiener processes are sampled on fine temporal grids, and the resulting Galerkin systems are advanced using problem-specific numerical integrators. At each time step, nonlinear terms are evaluated on spatial quadrature grids and projected back onto the retained spectral modes. The resulting spectral coefficients are then reconstructed at the evaluation times to obtain the reference trajectories $u_{\mathrm{ref}}$. All methods are evaluated along the same underlying Wiener paths. Reference accuracy is assessed by comparing time resolutions under the same noise realizations at fixed spatial resolution.

\section{Computational Cost} \label{app:com_cost}
In this section, we further analyze the computational cost of BRiSE across the six SPDE benchmarks. We separately report the offline training cost and the test-time solution cost. The former includes the complete training budget and periodic physics-based validation, while the latter measures the time required to solve a batch of eight held-out Wiener paths. Unlike BRiSE, which performs forward inference with a pretrained model, PINN and NeS-PINN require path-specific optimization at test time. We therefore report their optimization costs separately to distinguish them from forward-solve latency.

\subsection{Offline Training Cost of BRiSE}
We report the offline training cost of BRiSE across six SPDE benchmarks in Table~\ref{tab:training_cost}. All models are trained for 1,500 optimization steps on a single NVIDIA GeForce RTX 4090 GPU. The reported wall-clock training time includes the complete optimization process and periodic physics-based validation. As shown in Table~\ref{tab:training_cost}, training takes between 65.6 and 352 seconds across the six benchmarks, with model sizes ranging from approximately 41K to 1.55M parameters. The Allen--Cahn benchmark incurs the highest training cost due to its finer temporal discretization ($N=100$), compared with $N=20$ for the remaining benchmarks.

\begin{table}[htbp]
\centering
\caption{Offline training cost of BRiSE across six SPDE benchmarks.}
\label{tab:training_cost}
\begin{tabular}{lrrrr}
\toprule
\textbf{SPDE} & \textbf{Parameters} & \textbf{Intervals $\boldsymbol{N}$} & \textbf{Iters} & \textbf{Training (s)} \\
\midrule
Burgers        & 41,408    & 20  & 1,500 & 65.6 \\
Heat           & 41,408    & 20  & 1,500 & 70.7 \\
Fisher--KPP    & 42,178    & 20  & 1,500 & 70.0 \\
Allen--Cahn    & 41,408    & 100 & 1,500 & 352  \\
Cahn--Hilliard & 90,944    & 20  & 1,500 & 82.0 \\
Gray--Scott & 1,546,244 & 20  & 1,500 & 134  \\
\bottomrule
\end{tabular}
\end{table}

\subsection{Computational Cost at Test Time}
We compare the test-time computational costs of BRiSE with five baselines across six SPDE benchmarks. BRiSE performs forward inference using a pretrained model, whereas classical numerical solvers directly integrate each stochastic trajectory. In contrast, PINN and NeS-PINN require path-specific optimization at test time. We therefore distinguish forward-solve times from optimization-based solution times in Table~\ref{tab:test_time_cost}. BRiSE achieves substantially lower test-time costs than the two physics-informed neural baselines, although classical numerical solvers generally remain faster per solve. These results highlight the advantage of BRiSE in efficiently generating new stochastic trajectories without repeated neural network optimization.

\begin{table}[htbp]
\centering
\caption{Computational cost at test time across six SPDE benchmarks,
reported in milliseconds for a batch of eight Wiener paths.}
\label{tab:test_time_cost}
\begin{tabular}{lrrrrrr}
\toprule
Method
& Burgers
& Heat
& Fisher--KPP
& Allen--Cahn
& Cahn--Hilliard
& Gray--Scott \\
\midrule
\multicolumn{7}{l}{\textit{Pretrained inference / Numerical solvers}} \\
BRiSE
& 8.23 & 7.95 & 30.2 & 39.9 & 11.2 & 17.0 \\
EM
& 1.95 & 1210 & 3.50 & 9.27 & 2.93 & 6.41 \\
ExpEuler
& 2.00 & 2.11 & 2.45 & 9.52 & 2.37 & 6.55 \\
IMEX Euler
& 1.95 & 2.13 & 2.67 & 9.37 & 2.92 & 6.51 \\
\midrule
\multicolumn{7}{l}{\textit{Path-specific optimization}} \\
PINN
& 166000 & 83200 & 86200 & 306000 & 93200 & 1550000 \\
NeS-PINN
& 10100 & 10100 & 11000 & 48100 & 13800 & 79000 \\
\bottomrule
\end{tabular}
\end{table}

\end{document}